\documentclass[onefignum,onetabnum]{siamart171218}

\usepackage{lipsum}
\usepackage{amsfonts}
\usepackage{graphicx}
\usepackage{epstopdf}
\usepackage{algorithmic}
\usepackage{xcolor}
\usepackage{hyperref}
\ifpdf
  \DeclareGraphicsExtensions{.eps,.pdf,.png,.jpg}
\else
  \DeclareGraphicsExtensions{.eps}
\fi

\newcommand{\blue}[1]{{\color{black} #1}}
\newcommand{\green}[1]{{\color{black} #1}}
\newcommand{\R}{\mathbb{R}}
\newcommand{\dt}{{\Delta t}}
\newcommand{\du}{\Delta U}
\newcommand{\dw}{\Delta W}

\newcommand{\dV}{\Delta V}
\newcommand{\A}{\mathcal{A}}
\newcommand{\X}{\mathcal{X}}

\newtheorem{thm}{Theorem}[section]

\newtheorem{lem}[thm]{Lemma}

\newsiamremark{remark}{Remark}
\newsiamremark{hypothesis}{Hypothesis}
\crefname{hypothesis}{Hypothesis}{Hypotheses}
\newsiamthm{claim}{Claim}

\headers{Strong convergence analysis  for Langevin Dynamics }{X. Li and A. Telatovich}

\title{Strong Convergence Analysis for Numerical Methods for the Langevin Dynamics Model and Higher-order Operator-splitting Schemes
\thanks{Submitted to the editors DATE.
\funding{The research of X. Li is funded by NSF grant NSF Grant DMS-1522617 and DMS-1619661.}}}

\author{Xiantao Li\thanks{
  (\email{xli@math.psu.edu}, \url{http://personal.psu.edu/xxl12/}).}
\and Adam Telatovich\footnotemark[2]\thanks{
  (\email{atelatov@umd.edu}, \email{adam.telatovich@gmail.com}).}}

\usepackage{amsopn}

\ifpdf
\hypersetup{
  pdftitle={Strong Convergence Analysis for Numerical Methods for the Langevin Dynamics Model and Higher-order Operator-splitting Schemes},
  pdfauthor={X. Li and A. Telatovich}
}
\fi

\begin{document}

\maketitle

\begin{abstract}
  We provide the analysis of the strong convergence of some widely implemented methods for the Langevin dynamics model with additive noise. We show that a direct splitting of deterministic  and random terms, including the symmetric splitting methods, only offers strong convergence of order 1.  We present some new methods with higher order strong convergence. The new class of methods are operator-splitting schemes,  constructed based on the Kunita's solution representation for stochastic differential equations. We present stochastic algorithms with strong orders up to 3. Both mathematical analysis and numerical evidence are provided to verify the desired order of accuracy. We also examine the over-damped and under-damped limits of some of the new methods.
\end{abstract}



\section{Introduction}
The Langevin dynamics (LD) equation plays a fundamental role in the modeling of many complex dynamical  systems subject to random noise. In its simplest form, it can be expressed as the Newton's equations of motion with added frictional and random forces, which are usually introduced to model the influence of the surrounding environment, and posed to satisfy the fluctuation-dissipation theorem.


It does not come as a surprise that such equations only have explicit solutions in very rare cases. In general, approximate solutions have to be constructed at discrete time steps.  As a system of stochastic differential equations (SDE), there are various classical methods for approximating the solutions \cite{Kloeden1991}. However, low order methods, such as the Euler-Maruyama method, often do not have sufficient accuracy for accurate predictions. On the other hand, higher order methods that are constructed based on direct expansions of solutions (It\^o-Taylor expansions) usually involve  high order derivatives of the drift and diffusion coefficients, which makes the implementation rather difficult.   For instance, for bio-molecular models \cite{schlick2010molecular}, this implies that one has to compute the derivatives of the inter-molecular forces, which typically is not plausible. Extensions of Runge-Kutta methods, due to the many evaluations of the drift and diffusion terms at intermediate stages, have been largely neglected in the molecular simulation community. In molecular models \cite[Chapter 7]{leimkuhler2015molecular}, instead of using Runge-Kutta methods \cite{Rumelin1982}, operator-splitting schemes have been more widely implemented.
Such algorithms, especially with applications to molecular dynamics simulations,
have been treated extensively in \cite{leimkuhler2002molecular,leimkuhler2015molecular}, where many theoretical and practical aspects have been discussed.
The idea is to separate out terms on the right hand side and form two  or more SDEs, each of which can be solved explicitly. This is denoted by an [ABO] notation in  \cite{leimkuhler2002molecular}. Many existing methods can be recast into this form \cite{leimkuhler2013rational,leimkuhler2013robust,bou2010long,bussi2007accurate,melchionna2007design,brunger1984stochastic}.  One particular advantage of the splitting methods is that   they are very easy to implement, since each substep can be carried out exactly. The splitting methods can also be designed to better sample the equilibrium averages.   Another important  approach is based on solving the coordinate and momentum equations consecutively. For example, one can start by assuming the coordinates remain constant, and integrate out the momentum equation exactly. Then using this solution for the momentum, one can integrate the first equation and obtain an updated  coordinate for the next step. These two steps can be written in an operator splitting form. But a further correction can be made by assuming the force is linear in time, constructed using the coordinates at the current and next steps. This led to the stochastic velocity Verlet method (SVV) \cite{dellago1998efficient,allen1989computer,van1982algorithms}, which has been implemented in simulation packages, e.g., TINKER \cite{ponder2004tinker}. Other integration methods can also be found in the literature \cite{vanden2006second,skeel1999integration,ma2016fluctuation,melchionna2007design,nakajima1996new}.

On the theoretical side,  the fundamental issue of the numerical accuracy has been thoroughly discussed in \cite{leimkuhler2002molecular,cances2007theoretical}. In particular, the weak convergence of the numerical solution has been rigorously proved in \cite{leimkuhler2016computation}. Such analysis is crucial when the approximation methods are used to sample the corresponding {\it equilibrium} statistics. This is particularly useful when the averages of certain quantities are of interest. On the other hand, to the best of our knowledge, the theoretical analysis of {\it strong } convergence has not been fully studied. Strong convergence ensures the accuracy in terms of individual realizations and solutions at transient stages \cite{sauer2012numerical}. Such notation of accuracy would be useful in non-equilibrium simulations, such as non-equilibrium molecular dynamics (NEMD) simulations, which have been very useful in the study of transport processes \cite{im2002ion,mashreghi2011molecular,isralewitz2001steered,ilg2005structure}.
Strong convergence usually implies weak convergence, but not vice versa.
Typically, strong convergence can be examined by comparing to the It\^{o}-Taylor expansion. Therefore, the fact that the splitting methods discussed in the literature often do not involve multiple It\^o integrals of order 2 or higher is already an indication that those methods are only of strong order 1 or less, regardless of how the operators are split up and how many fractional steps are involved.

This paper will primarily focus on this theoretical aspect. In the first half, we analyze two widely used methods, namely the SVV and a naive splitting method. In the second half, we present several new operator splitting methods with higher {\it strong} order.
Our starting point for the new methods is the solution representation by Kunita \cite{kunita1980}. Written formally as an operator exponential form, the differential operator is expressed in terms of the commutators involving the differential operators associated with the drift and diffusion coefficients, along with multiple It\^o integrals. Intuitively, we can make truncations at various levels, yielding approximation methods of increasing order. Such truncation schemes have been used in \cite{Misawa2000} for one-dimensional stochastic differential equations as a starting point to construct robust algorithms for scalar SDEs with multiplicative noise. It was demonstrated that such algorithms can preserve the non-negativity of the solution. In this paper, we extend the applications of these truncations to the Langevin dynamics, introducing truncations of the solution operator, and obtaining approximate solutions that can be written as solutions of ODEs, for which many efficient methods exist. We choose the well established operator splitting methods for these ODEs, yielding new operator splitting schemes for the Langevin dynamics model.  With direct mathematical analysis and numerical tests, we verify that  these methods have strong order 2 and 3. They can be used as highly accurate numerical tools to study stochastic dynamics modeled by Langevin dynamics.
No derivatives of the force function are required for the order 2 splitting method,
and only first order derivatives are needed for the order 3 splitting method.


The rest of the paper is organized as follows. \cref{sec: basics} presents the basic theory for strong convergence. In \cref{sec: analysis}, we analyze the convergence of some existing numerical methods. In \cref{sec: os}, we introduce the new class of operator splitting schemes based on the truncations of the Kunita's solution operator and examined the strong order of accuracy.
\cref{sec: tests} contains numerical tests that will demonstrate the expected order of convergence.
\cref{sec: proofs} contains our proofs of the convergence results we presented in \cref{sec: analysis} and \cref{sec: os}, and \cref{sec: IT} is the Appendix \green{where we derive the It\^o-Taylor expansions for the Langevin model}.

\section{The basic theory}\label{sec: basics}
          This paper is mainly concerned with the Langevin dynamics model with $n$ spacial dimensions,
\begin{equation}\label{eq: lgv}
 \left\{\begin{aligned}
   dx=& vdt,\\
   dv=& f(x)dt -\Gamma vdt + \sigma dW_t,
\end{aligned}\right.
\end{equation}
where $x=(x^1,\dots,x^n), v=(v^1,\dots,v^n) \in \mathbb{R}^n$ can be interpreted as position and velocity components respectively,  $W(t)=(W^1(t),\dots,W^n(t))\in\R^n$ is the standard $n-$dimensional Brownian motion and $0\leq t\leq T$. $W'(t)$ represents white noise in time. Assume the function $f=f(x):\R^n\rightarrow\R^n$, representing the conservative force (for example, the Morse potential), has bounded second derivatives. Here we will consider the case where $\sigma$ and $\Gamma$ are constant $n \times n$ matrices which we assume satisfy the fluctuation dissipation theorem,
\begin{equation}\label{eq: FDT}
    \sigma\sigma^T = 2\Gamma k_B T,
\end{equation}
where $k_B$ is Boltzman's constant, $\sigma^T$ denotes the transpose of $\sigma$, and $T$ represents the temperature of the system. In particular, the noise is additive, and as a consequence, \eqref{eq: lgv} may be interpreted as either an It\^o or Stratonovich SDE. For particle dynamics, typically $n=3N$, with $N$ being the total number of particles.

We can view the Langevin equation above as an example of an autonomous It\^o SDE
\begin{equation}\begin{aligned}
    dz_t = a(z_t)dt + bdW_t,
\end{aligned}\end{equation}
where $z_t=(x_t,v_t)$, and $a,b$ are the respective drift and (constant) diffusion terms, given here by
\begin{equation}\begin{aligned}
    a=\left(\begin{matrix}v\\ f(x)-\Gamma v\end{matrix}\right)\in \R^{2n},
        \text{\quad}b = \left(\begin{matrix} 0\\ \sigma \end{matrix}\right)\in \R^{2n\times n}.
\end{aligned}\end{equation}

The numerical solutions will involve Brownian increments, whose $j$th component ($j=1,\dots,n$) at the $k$th step ($k=0,\dots,n_T-1$) are $\sqrt{\dt}N(0,1)$ random variables $\dw_k^j:= W_{t_{k+1}}^j - W_{t_{k}}^j$, where $\{t_k\}_{k=1}^{n_T}$ is the uniform discretization $t_k=k\dt$ of a time interval $[0,T]$ with uniform step size $\dt = t_{k+1}-t_k$. In vector form, we write $\dw_k = (\dw_k^1,\dots,\dw_k^n)$.
The components $\dw_k^i$ and $\dw_k^j$ are pairwise independent for each step $k$, and the increments $\dw_k^j$ and $\dw_l^j$ are independent for each component $j$.

A discrete time approximation $Y=Y^{\dt}$ of $X$ with uniform step size $\dt$ converges with strong order $\gamma\in\{0.5, 1, 1.5, 2,\dots\}$ at time $T>0$ provided that  there are  constants $C$ and $\Delta>0$ such that
\begin{equation}
\mathbb{E}\left(\sup_{0\leq k\leq n_T}|X(k\dt)-Y(k\dt)|\right)\leq C \dt^{\gamma}\text{\quad for all }0<\dt<\Delta.
\end{equation}
The notation $\mathbb{E}(\cdot)$ denotes the usual expectation of a random variable on a probability space $(\Omega,\mathcal{F},P)$. The random variables $X,Y$ are driven by the same process $W_t$.
Clearly, strong convergence is related to the convergence of paths.
\begin{remark}
We point out that the definition above actually means that the convergence is uniform over the entire time interval $[0,T]$,
whereas the definition given in \cite[eq. 9.6.3]{Kloeden1991} means convergence at the terminal time $T$:
\begin{equation}
    \mathbb{E}\left( |X(T)-Y(T)| \right)\leq C\dt^\gamma \text{\quad for all }0<\dt<\Delta.
\end{equation}
We will use the former definition of strong convergence instead of the latter, showing preference for uniform
estimates over point-wise estimates.
\end{remark}


\section{\green{Analysis of some existing methods}}\label{sec: analysis}

\subsection{It\^o-Taylor expansion of the solution}
To construct an algorithm with certain strong order of accuracy,  a direct comparison has to be made with
the It\^{o}-Taylor expansion of the exact solution  \cite{Kloeden1991}. With direct, but involved calculations, we have
obtained such expansions for the  Langevin dynamics model \eqref{eq: lgv}, given as follows,
\begin{lem}
The strong It\^{o}-Taylor approximations of order 1, 2 and 3 for the Langevin equation \green{\eqref{eq: lgv}} with additive noise are one-step methods given at step $k+1$, for $k=0,1,\dots,n_T-1$, respectively, by
\begin{equation}\label{eq: IT0}
\begin{aligned}
    \left(\begin{matrix} x_{k+1}\\v_{k+1}\end{matrix}\right)&=\left( \begin{matrix} x_k\\v_k \end{matrix}\right)+\left( \begin{matrix} v_k\\f(x_k)-\Gamma v_k \end{matrix} \right)\dt+  \left( \begin{matrix} 0\\\sigma\dw_k \end{matrix} \right) + H.O.T.,\\ \\
    \left(\begin{matrix} x_{k+1}\\v_{k+1}\end{matrix}\right)&=\left( \begin{matrix} x_k\\v_k \end{matrix}\right)+\left( \begin{matrix} v_k\\f(x_k)-\Gamma v_k \end{matrix} \right)\dt
     + \left( \begin{matrix} 0\\\sigma\dw_k \end{matrix} \right)\\ \\
     &   +\left( \begin{matrix} f(x_k)-\Gamma v_k\\  \nabla_x \left(f(x_k)\right)v_k-\Gamma f(x_k)+\Gamma^2 v_k \end{matrix} \right)\frac{\dt^2}{2}
        + \left( \begin{matrix} \sigma I_{(J,0),k}(\dt)\\-\Gamma\sigma I_{(J,0),k}(\dt) \end{matrix} \right)  + H.O.T.,\\ \\
    \left(\begin{matrix} x_{k+1}\\v_{k+1}\end{matrix}\right)&=\left( \begin{matrix} x_k\\v_k \end{matrix}\right)+\left( \begin{matrix} v_k\\f(x_k)-\Gamma v_k \end{matrix} \right)\dt
    + \left( \begin{matrix} 0\\\sigma\dw_k \end{matrix} \right)\\ \\
   & +\left( \begin{matrix} f(x_k)-\Gamma v_k\\  \nabla_x  \left(f(x_k)\right)v_k-\Gamma f(x_k)+\Gamma^2 v_k \end{matrix} \right)\frac{\dt^2}{2}
    + \left( \begin{matrix} \sigma I_{(J,0),k}(\dt) \\-\Gamma\sigma I_{(J,0),k}(\dt) \end{matrix} \right)\\ \\
    &+\left(\begin{matrix} \nabla_x\left(f(x_k)\right)v_k-\Gamma f(x_k)+\Gamma^2 v_k\\
                 \nabla_x[\nabla_x(f(x_k))v_k]v_k- \Gamma \nabla_x(f(x_k))v_k \end{matrix}\right)\frac{\dt^3}{3!}\\ \\
    &+\left(\begin{matrix} 0 \\  \nabla_x(f(x_k))(f(x_k)-\Gamma v_k) +\Gamma^2(f(x_k)-\Gamma v_k)  \end{matrix}\right)\frac{\dt^3}{3!} \\ \\
    &+\left(\begin{matrix} -\Gamma\sigma I_{(J,0,0),k}(\dt)\\ \left( \nabla_x(f(x_k))\sigma +\Gamma^2\sigma \right) I_{(J,0,0),k}(\dt) \end{matrix}\right)+ H.O.T.
\end{aligned}\end{equation}
where $ \sigma^j$ is the $j$th column of $\sigma\in\R^{n\times m}$, $ \nabla_x  f(x)\in\R^{n\times n}$ is the
gradient matrix of $f(x)\in\R^n$ at $x$, $(\nabla_x f(x))_{ij}=\frac{\partial f_i(x)}{\partial x^j}$, and
$\nabla_x[ \nabla_x (f(x))v ]$ is the gradient matrix of $\nabla_x f(x)v\in\R^n$ at $x$.
The abbreviation H.O.T. stands for Higher Order Terms (with respect to $\dt$).
For the calculations involved to derive the expansions above, the reader may refer to \cref{sec: IT}.
Furthermore, the strong order 0.5, 1.5 and 2.5 It\^{o}-Taylor approximations are the same as the strong order 1, 2 and 3 approximations, respectively, due to the fact that the noise is additive.
\end{lem}

\green{It is important to point out that} $I_{(J,0),k}=(I_{(1,0),k},\dots,I_{(n,0),k})\in\R^n$ and the $I_{(J,0,0),k}=(I_{(1,0,0),k},\dots,I_{(n,0,0),k})\in\R^n$ are random vectors, with components $I_{(j,0),k}$ and $I_{(j,0,0),k}$ which are multiple stochastic integrals, given by the following formulas
\begin{equation}\begin{aligned}
    I_{(j,0),k}(t)&=\int_{t_k}^{t_{k+1}} W_s^j ds \\ I_{(j,0,0),k}(t)&= \int_{t_k}^{t_{k+1}}\int_{t_k}^{s}W_u^j du ds.
\end{aligned}\end{equation}
The \green{quantities} $\dw_k$, $I_{(J,0),k}$ and $I_{(J,0,0),k}$ are independent Gaussian
random variables, whose means and variances do not depend on $k$ due to translation invariance, and they
have the same distributions, respectively, as the following random variables:
\begin{equation}
    \dw := W_{\dt},\text{\quad}I_{(j,0)}:=\int_0^{\dt}W_s^j ds,\text{\quad and\quad}I_{(j,0,0)}:=\int_0^{\dt}\int_0^s W_u^j du ds.
\end{equation}
The expansions \green{ \eqref{eq: IT0} } reduces to the usual Taylor expansions when $\sigma = 0$. We refer the readers
to \cref{sec: IT} or to \cite{Kloeden1991} for the detailed explanation of the notations. How the terms in the expansion are retained is determined by an estimate on the variance of each term, and how they accumulate in time \cite{Kloeden1991}. The additive noise has the desirable effect that several terms vanish in the It\^o Taylor expansion.

One-step numerical schemes can be directly obtained from the It\^{o}-Taylor expansion \eqref{eq: IT0}.
For example, by neglecting the high order terms (H.O.T.) in the first expansion, one obtains the
familiar Euler-Maruyama (EM) method, which in this case has strong order 1. It is tempting to make
truncations for the order 2 and 3 expansions. However, those direct expansions involve higher order
derivatives of $f(x)$, and they may not be easy to implement. Specifically, the order 2 expansion
requires computing $\nabla_x f$ and the order 3 expansion requires both $\nabla_x f$ and $\nabla_x^2 f$.
Therefore, we only use these expansion as a guideline to \green{prove strong order convergence and} construct algorithms with high strong order,
in light of the theoretical analysis (Theorem 11.5.1 \cite{Kloeden1991}).

%
%

\subsection{\green{Some existing splitting methods}}

A natural approximation of \eqref{eq: lgv} can be obtained by splitting the equation into several subproblems, each of which can be solved exactly. A wide variety of splitting methods have been discussed in \cite{leimkuhler2002molecular} \green{that have this flavor}. For example, we may consider to  split the Langevin equation as follows,
\begin{equation}\label{eq: A}\left\{
\begin{aligned}
   x'=& v\\
   v'=& 0,
\end{aligned}\right.
\end{equation}
and
\begin{equation}\label{eq: B}\left\{
 \begin{aligned}
   x'=& 0\\
   v'=& f(x) -\Gamma v + \sigma W'(t).
\end{aligned}\right.
\end{equation}
Both of these equations have explicit solutions. The solution steps can be denoted by abstract operators $A$ and $B$, respectively, and approximations can be obtained by following the operations, e.g., $AB$, $\frac A2 B \frac A2$, etc \cite{leimkuhler2002molecular}. Unfortunately, the strong accuracy of these methods is very limited, \green{as shown by the following theorem}.

\begin{thm}
The symmetric and non-symmetric splitting methods,
$$\exp(A/2)\exp(B)\exp(A/2) \text{\quad and\quad} \exp(A)\exp(B),$$
respectively, have strong order 1.
\end{thm}
The proof will be postponed to \cref{sec: pfdirectsplit}.

\bigskip
\subsection{\green{The stochastic velocity-Verlet (SVV) method}}

Another widely implemented scheme is the stochastic velocity-Verlet's (SVV) method \cite{van1982algorithms},
as we mentioned in the introduction. This method starts with an assumption that $x(t)$ remains as a constant,
and integrates the second equation in \eqref{eq: lgv} exactly, giving rise to,\green{
\begin{equation}
 v(t) = c_0(t)v_k + c_1(t)f(x_k) + \int_{t_k}^t e^{-\Gamma (t-s)} \sigma dW_s,
\end{equation}}
where $t_k\leq t\leq t_{k+1}$, \green{the time-dependent coefficients $c_0,c_1$ are given by,
\begin{equation}\begin{aligned}
    c_0(t)&= e^{-\Gamma t},\\
    c_1(t) &= \Gamma^{-1}(I- c_0(t)),
\end{aligned}\end{equation}}
$I$ denotes the $n\times n$ identity matrix and $e^{-\Gamma t}$ is a
matrix exponential. (For a derivation of this formula, see \cref{sec: pfdirectsplit}.)

With this approximation of $v(t)$, one can now turn to the first equation, and integrate. This gives,
\begin{equation}\label{eq: svv-x}
 x_{k+1} = x_k + c_1(\dt) v_k + c_2(\dt) f(x_k) + \int_{t_k}^{t_{k+1}}   \int_{t_k}^t e^{-\Gamma (t-s)} \sigma dW_s dt.
\end{equation}
Here the new coefficient $c_2$ is given by,\green{
\begin{equation}
    c_2(\dt)=  \int_{t_k}^{t_k+\dt}  \Gamma^{-1}(I- e^{-\Gamma s})ds.
\end{equation}}
One might stop here and accept the position and velocity values. But in SVV one  uses the updated position value and approximate the function $f$
 by a linear function,
 \begin{equation}
 f(x(t)) \approx f(x(0)) + \frac{\big(f(x(\dt))-f(x(0))\big) t}{\green{\dt}}.
\end{equation}
With this approximation, one can integrate the velocity equation again. One finds that,
\begin{equation}\label{eq: svv-v}
 v_{k+1} = c_0\green{(\dt)} v_k + c_1\green{(\dt)} f(x_k) + c_2\green{(\dt)} \big(f(x_{k+1})-f(x_k)\big)+ \int_0^\dt e^{-\Gamma (\dt-s)} \sigma dW_s.
\end{equation}

Equations \eqref{eq: svv-x} and \eqref{eq: svv-v} form the basis for the SVV method.
The formulas can be repeated, and at each step, the function $f$ is evaluated only once at each step,
which is typically viewed as a considerable advantage.
As we will prove, this algorithm has second order strong accuracy.

\begin{thm}
  The SVV algorithm has strong order 2.
\end{thm}
A brief proof can be found in \cref{sec: svv}. To the best of our knowledge, this is a new theoretical result.


\section{New operator-splitting algorithm with higher order strong convergence}\label{sec: os}


Here we propose new splitting algorithms. Our starting point is the Kunita's solution operator \cite{kunita1980}. In particular, for the standard SDE written in the differential form,
\begin{equation}\label{eq: SDEs}
  d z_t = a(z_t) dt + b(z_t) dW_t,\text{\quad} z(0)=z,
\end{equation}
we define the differential operators,
\begin{equation}
 \mathcal{X}_0 = a \cdot \nabla_{z},  \quad \mathcal{X}_j= b^j \cdot \nabla_{z},
\end{equation}
Here $z_t $ represents the solution at time $t$ and $b^j$ denotes the $j$th column of $b$. Then the exact solution of the SDE can be formally expressed as,
\begin{equation}
 z_t = \exp(D_{t})z,
\end{equation}
where, by \cite[eq. (2.5)]{Misawa2000}, for $t=\dt$ we have,
\begin{equation}\begin{aligned}
 D_{\dt} &=\dt \mathcal{X}_0 +  \sum_{j=1}^n \dw^j  \mathcal{X}_j + \frac{1}{2}\sum_{j=1}^n[\dt,\dw^j][ \mathcal{X}_0, \mathcal{X}_j] \\
  &+ \frac{1}{18}\sum_{j=1}^n [[\dt,\dw^j],\dt][[\mathcal{X}_0,\mathcal{X}_j],\mathcal{X}_0]+\dots
\end{aligned}\end{equation}
Here, $[\mathcal{X}_0,\mathcal{X}_j]$ denotes the commutator bracket of differential operators:
\begin{equation}
\mathcal{X}_0 \mathcal{X}_j-\mathcal{X}_j \mathcal{X}_0,
\end{equation}
and $[[\mathcal{X}_0,\mathcal{X}_j],\mathcal{X}_0]=[\mathcal{X}_0,\mathcal{X}_j]\mathcal{X}_0-\mathcal{X}_0 [\mathcal{X}_0,\mathcal{X}_j]$.
Higher order commutators can be defined similarly.

The \green{corresponding terms} $[\dt,\dw^j]$ \green{(for $j=1,\dots,n$)} are \green{combinations of double stochastic integrals} given by
\begin{equation}
[\dt,\dw_k^j]:=\int_{t_k}^{t_{k+1}} tdW_t^j-\int_{t_k}^{t_{k+1}} W_t^j dt,
\end{equation}
and the \green{corresponding terms} $[[\dt,\dw_k^j],\dt]$ \green{(for $j=1,\dots,n$)} are \green{combinations of triple stochastic integrals} given by, \cite[eq. (2.15)]{Misawa2000}
\begin{equation}
    [[\dt,\dw_k^j],\dt]=\frac{1}{18}(2I_{(0,j,0),k}-2I_{(j,0,0),k}+\dt I_{(j,0),k}-\dt I_{(0,j),k}).
\end{equation}
\green{The definitions of $I_{(j,0)}$ and $I_{(j,0,0)}$ are double and triple stochastic integrals, defined in \cref{sec: analysis} for the It\^o-Taylor expansions \eqref{eq: IT0}.} The $I_{(0,j)}$ and $I_{(0,j,0)}$ \green{(which did not show up in the It\^o Taylor expansions)} are defined similarly:
\begin{equation}
    I_{(0,j),k}=\int_{t_k}^{t_{k+1}}\int_{t_k}^{s}dudW_s^j,\text{\quad}
        I_{(0,j,0),k}=\int_{t_k}^{t_{k+1}}\int_{t_k}^{s}\int_{t_k}^u dvdW_u^j ds,
\end{equation}
for $j=1,\dots,m$. These integrals do not show up in the It\^o-Taylor expansions \eqref{eq: IT0}
because their corresponding coefficient functions $g_{(0,j)}$ and $g_{(0,j,0)}$
(defined in \cite[Sec. 5.3]{Kloeden1991}) are identically zero for the Langevin dynamics with additive noise.
Indeed, one can make the  following calculation,
\begin{equation}
    g_{(0,j)}(x,v)=\left(\begin{matrix}0_{n\times 1} \\ \nabla_x(\sigma^j)v+\nabla_v(\sigma^j)(f(x)-\gamma v)\end{matrix}\right)=0_{2n\times 1}
\end{equation}
since the $\sigma^j$ are constant, and
\begin{equation}
    g_{(0,j,0)}(x,v)=\left(\begin{matrix} \nabla_x(\sigma^j)v+\nabla_v(\sigma^j)(f(x)-\gamma v)\\
            \nabla_x(-\gamma\sigma^j)v+\nabla_v(-\gamma\sigma^j)(f(x)-\gamma v)\end{matrix}\right)
            =0_{2n\times 1}
\end{equation}
since $\gamma\sigma^j$ is constant as well. The third order commutators above are used in \cref{sec: trunc3}.

For implementation, we use the fact that the stochastic integrals $[\dt,\dw_k^j]$ are independent and identically distributed (for all $k$), with the same distribution as
\begin{equation}
    \int_0^{\dt}tdW_t^j - \int_0^{\dt} W_t^j dt.
\end{equation}

Here we have used $\exp(D_{\dt})$ to define the numerical solution at $t=\dt$.

\green{In addition to the analysis, we also performed numerical tests, including a pendulum model and a Lennard-Jones cluster. The details will be described in \cref{sec: tests}.}


\subsection{First-order truncation}
We first make a truncation and keep the first two terms \cite[eq. 3.18]{Misawa2000}:
\begin{equation}\label{eq: D1}
 \begin{aligned}
D_{\dt}^I &= \dt \mathcal{X}_0 + \sum_{j=1}^n \dw^j \mathcal{X}_j\\
&= \dt v\cdot \nabla_x + \left( \dt(f(x)-\Gamma v)+\sum_{j=1}^n \sigma^j \dw^j \right)\cdot\nabla_v.
\end{aligned}
\end{equation}
One step of the numerical solution consists in integrating \eqref{eq: D1} over time $0\leq t\leq 1$.

 Once the Brownian motion $\dw$ has been sampled (and realized), the operator $\exp(D_{\dt}^I)$ corresponds to the
 solution operator of the following ODE system,
 \begin{equation}\label{eq: ODE-I}
	\begin{cases}
		x' &= \dt v\\
		v' &= \dt( f(x)-\Gamma v ) + \sum_{j=1}^n \sigma^j \dw^j,
	\end{cases}
\end{equation}
which we solve over the interval $0\leq t\leq 1$.
This approximation by the solution of the above ODE system will be referred to as {\bf truncation I}.

At this point, we can prove the strong order convergence \green{of the approximation using \eqref{eq: ODE-I}}. To see the local consistency, we expand the solution of the ODEs at $t=0,$
\begin{equation}\begin{aligned}
\exp(D_\dt^I)\left( \begin{matrix} x_k\\ v_k \end{matrix} \right)&=
    \left( \begin{matrix} x_k\\ v_k \end{matrix} \right)+\left( \begin{matrix} v_k\\f(x_k)-\Gamma v_k \end{matrix} \right)\dt+\sum_{j=1}^n\left( \begin{matrix} 0\\ \sigma^j \end{matrix} \right) \dw^j\\
    &+\left( \begin{matrix} f(x_k)-\Gamma v_k\\ \nabla_x f(x_k)v_k-\Gamma(f(x_k)-\Gamma v_k) \end{matrix} \right)\frac{\dt^2}{2}\\
    &+\sum_{j=1}^n \left( \begin{matrix} \sigma^j\\-\Gamma\sigma^j \end{matrix} \right) \frac{\dt\dw^j}{2}+H.O.T.
\end{aligned}\end{equation}
where the higher order terms do not involve $I_{(j,0)}$. (See \cref{sec: pftrunc1} for a derivation of this expansion.)
Thus with a comparison to \eqref{eq: IT0} we have,

\begin{thm}
For the Langevin equation with additive noise, the truncation method given by $z_{k+1}=\exp(D_\dt^I)z_k$ is precisely a strong order 1   approximation.
\end{thm}
See \cref{sec: pftrunc1} for the proof.

To solve the ODEs, we consider $D_{\dt}^{\rm I}=A+B$ where
\begin{equation}
    A=\dt v\cdot\nabla_x\text{\quad and\quad} B=\left((f(x)-\Gamma v)\dt+\sigma\dw \right)\cdot\nabla_v.
\end{equation}
By the Baker Campbell Hausdorff (henceforth BCH) formulas \cite[eq. 3.1]{Yoshida1990},
along with a direct comparison with the It\^{o}-Taylor expansion, we have established the following result,
\begin{thm}\label{thm: order1}
The non-symmetric splitting scheme
    \begin{equation}
        \exp \left(D^{\rm I}\right)\approx\exp (A)\exp (B)
    \end{equation}
and symmetric splitting scheme
    \begin{equation}
        \exp\left(D^{\rm I}\right)\approx\exp(A/2)\exp(B)\exp(A/2)
    \end{equation}
both yield approximations with strong order $\gamma=1$.
\end{thm}
\blue{\begin{proof}
    We give a short proof in \cref{sec: pftrunc1} that $\exp D_{\dt}^{\rm I}$ converges
    to the exact solution $z_{k+1}=\exp D_{\dt}(z_k)$ with strong order $\gamma=1$, that is, locally, we have,
    \begin{equation}
        |\exp D_{\dt}(z_k)-\exp D_{\dt}^{\rm I}(z_k)|\leq C_1 \dt^2
    \end{equation}
    for some constant $C_1>0$. On the other hand, it is well known \cite{Hairer2005} that the non-symmetric
    splitting $\exp A\exp B(z_k)$ is a first order ODE approximation, that is,
    \begin{equation}
        |\exp D_{\dt}^{\rm I}(z_k)- \exp A\exp B(z_k)|\leq C_2\dt^2
    \end{equation}
    for some constant $C_2>0$. By the triangle inequality, the non-symmetric splitting yields the following
    local error estimate with the exact solution:
    \begin{equation}
    \begin{aligned}
        |\exp D_{\dt}(z_k)&-\exp(A)\exp(B)(z_k)|\\
        &\leq |\exp D_{\dt}(z_k)-\exp D_{\dt}^{\rm I}(z_k)|
            +|\exp D_{\dt}^{\rm I}(z_k)-\exp(A)\exp(B)(z_k)|\\
        &\leq C_1\dt^2 + C_2\dt^2= \tilde{C}_k\dt^2,
    \end{aligned}\end{equation}
    where $\tilde{C}_k=C_1+C_2$.
    Taking the supremum over all $0\leq k\leq n_T-1$ and adding the right-hand sides, we get,
    \begin{equation}\begin{aligned}
        \sup_{0\leq k<n_T}|\exp D_{\dt}(z_k)-\exp A\exp B(z_k)|&\leq \dt^2\sum_{k=0}^{n_T-1} \tilde{C}_k \\
            &\leq \dt^2\sum_{k=0}^{n_T-1} \max_k \tilde{C}_k\\
            &=\dt^2 n_T \max_{0\leq k<n_T} \tilde{C}_k = \max_{0\leq k<n_T} \tilde{C}_k\cdot \dt,
    \end{aligned}\end{equation}
    so that, taking expectations, we obtain,
    \begin{equation}
        \mathbb{E}\left( \sup_{0\leq k<n_T} |\exp D_{\dt}(z_k)-\exp A\exp B(z_k)| \right)\leq C\dt
    \end{equation}
    where $C:=\max_{0\leq k<n_T} \tilde{C}_k$.

    The symmetric splitting $\exp(A/2)\exp(B)\exp(A/2)$ is well-known as a second order operator splitting
    method for ODEs \cite{Hairer2005}, but due to the lack of the stochastic integrals $I_{(j,0)}$,
    it only yields a first order approximation.
\end{proof}}

In our numerical tests, we sample the increments $\dw_k=(\dw_k^1,\dots,\dw_k^n)$ at the $k$th step ($k=0,1,\dots,n_T-1$) as follows:
\begin{equation}
    \dw_k = \sqrt{\dt}\vec{\xi}_k,
\end{equation}
where $\vec{\xi}_k\in\R^n$ is a random vector in $N(0,I_{n\times n})$ and the increments $\vec{\xi}_k,\vec{\xi}_l$
are independent for $k\neq l$. From \cref{fig: pend} and \cref{fig: LJ},
we see that the (non-symmetric and symmetric) operator splitting methods applied to truncation $D_\dt^I$ both converge with order 1.
Therefore, to obtain higher order strong convergence, a further truncation is needed.


\subsection{Second-order truncation}
 Now we consider the truncation of $D_t$ which includes the first order bracket \cite[eq. 3.22]{Misawa2000}:
 \begin{equation}\label{eq: trunc2}
	D_{\dt}^{\rm II} = \dt \mathcal{X}_0 + \sum_{j=1}^n \dw^j \mathcal{X}_j + \frac{1}{2}\sum_{j=1}^n[\dt,\dw^j][\mathcal{X}_0,\mathcal{X}_j].
\end{equation}
\green{For the Langevin dynamics model \eqref{eq: lgv},} we can write the differential operator $[\mathcal{X}_0,\mathcal{X}_j]$ as
\begin{equation}\label{eq: bracket}
    [\mathcal{X}_0,\mathcal{X}_j]=\left( - \frac12\sum_{j=1}^n\sigma^j[\dt,\dw^j] \right)\cdot \nabla_x + \left( \frac12\sum_{j=1}^n \Gamma\sigma^j[\dt,\dw^j] \right)\cdot \nabla_v,
\end{equation}
which allows us to rewrite equation \eqref{eq: trunc2} above as,
\begin{equation}\begin{aligned}
	D_{\dt}^{\rm II}&= \left( \dt v - \frac12\sum_{j=1}^n\sigma^j[\dt,\dw^j] \right)\cdot \nabla_x \\
            &+ \left(\dt \big(f(x)-\Gamma v \big)+ \sum_{j=1}^n \sigma^j\dw^j + \frac12\sum_{j=1}^n \Gamma\sigma^j[\dt,\dw^j] \right)\cdot \nabla_v.
\end{aligned}\end{equation}
This will be referred to as {\bf truncation II}.

We note in passing that the expression \eqref{eq: bracket} for $[\mathcal{X}_0,\mathcal{X}_j]$ is valid
so long as the operator is applied to functions that are linear in $v$, which will always be the case in
this paper. And for functions linear in $v$, we note that $\mathcal{X}_0,\mathcal{X}_j$ equal the differential
operators $L^0,L^j$ respectively, which appear in the definition of It\^o-Taylor expansions
(see \cite[Sec. 5.1]{Kloeden1991}).

Now we define Gaussian random variables $\du =(\du^1,\dots,\du^n)$, where
\begin{equation}
\du^j:=  \frac{1}{2}[\dt,\dw^j] = I_{(0,j)}  -  \frac12 \dt \dw^j, \text{\quad}j=1,2,\cdots,n.
\end{equation}
\blue{The random variable $\du^j$ will play an important role in obtaining higher order strong convergence. For discretizations which lack higher order stochastic integrals (such as $I_{(0,j)}$ in the definition of $\du^j$), and only have the increments $\dw^j$, there is a barrier to higher order convergence. The Euler-Maruyama \green{method} obtains the highest convergence rate possible for such schemes, and this barrier is discussed in more detail in the paper by \cite[Thm. 3]{Rumelin1982}.  As we saw in the It\^o-Taylor schemes \eqref{eq: IT0}, the higher order schemes require higher order stochastic integrals.}

With $\dw:=(\dw^1,\dots,\dw^n)$, we have
\begin{equation}
D_{\dt}^{\rm II} = \left( \dt v -\sigma \du \right)\cdot \nabla_x + \left( \dt(f(x)-\Gamma v)+ \sigma\dw +\Gamma\sigma \du\right)\cdot\nabla_v
\end{equation}

Once $\dw$ and $\du$ are realized, the solution corresponds to that of the following ODEs at time $t=1$,
\begin{equation*}
	\begin{cases}
		x' &= v \dt -\sigma \Delta U,\\
		v' &= f(x) \dt -\Gamma v\dt +{\sigma\dw}+\Gamma\sigma\du.
	\end{cases}
\end{equation*}

A direct expansion of the solutions of the ODEs is given by,
    \begin{equation}\begin{aligned}
    \exp D^{\rm II}\left(\begin{matrix}x_k\\v_k \end{matrix}\right)
        &= \left(\begin{matrix}x_k\\v_k\end{matrix}\right)
            +\left(\begin{matrix}v_k\\ f(x_k)-\Gamma v_k \end{matrix}\right)\dt
            +\sum_{j=1}^n \left(\begin{matrix} 0_{n\times 1} \\ \sigma^j \end{matrix}\right)\dw_k^j\\
            &+\left(\begin{matrix}f(x_k)-\Gamma v_k\\ \nabla_x f(x_k)v_k-\Gamma f(x_k)+\Gamma^2 v_k \end{matrix}\right)\frac{\dt^2}{2}
            +\sum_{j=1}^n \left(\begin{matrix} \sigma^j \\ -\Gamma\sigma^j \end{matrix}\right)I_{(j,0),k}\\
            &+\sum_{j=1}^n \left(\begin{matrix} \Gamma\sigma^j \\ -\nabla_xf(x_k)\sigma^j-\Gamma^2\sigma^j \end{matrix}\right)\frac{\dt\du_k^j}{2}+\text{H.O.T.},
    \end{aligned}\end{equation}
    \blue{where H.O.T. stands for Higher Order Terms. We observe that the stochastic integrals $I_{(0,j)}$
    do not appear among the lower order terms. This is because the commutator $[\dt,\dw^j](=2\du^j)$
    in $D^{\rm II}$ allows us to include the stochastic integrals $I_{(j,0)}$ that we see in the
    higher order It\^o-Taylor expansions, while simultaneously excluding the unnecessary integrals $I_{(0,j)}$.
    A short derivation of the expansion above can be seen in \cref{sec: pftrunc2}.}

Thus we have,
\begin{thm}
 The operator $\exp\left(D_{\dt}^{\rm II}\right)$ generates a solution with strong order 2.
\end{thm}

For the numerical implementation, we use the following splitting, $D_{\dt}^{\rm II}=A+B$, where,
\begin{equation}
\begin{aligned}
  A=& \big( \dt v  -\sigma \Delta U\big) \cdot\nabla_x\\
  B= & \big( f(x) \dt -\Gamma v\dt +{\sigma\dw}+\sigma\Gamma\Delta U \big)\cdot \nabla_v.\\
\end{aligned}
\end{equation}

\blue{For practical implementations, we \green{denote} the increments by $\dw_k$ and $\du_k$ at the $k$th step ($k=0,1,\dots,n_T-1$), \green{and sample them} using the covariance matrix:
\begin{equation}\label{eq: covtrunc2}
C\overset{\text{def}}{=}     \left( \begin{matrix} \mathbb{E}(\dw_k^2)&\mathbb{E}(\dw_k\du_k)\\
                            \mathbb{E}(\du_k\dw_k)&\mathbb{E}(\du_k^2)\\
                            \end{matrix} \right)
    = \left(\begin{matrix} \dt I_{n\times n}& 0_{n\times n}\\ 
                            0_{n\times n}&\frac{\dt^3}{12}I_{n\times n} \\
                            \end{matrix}\right).
\end{equation}
This covariance formula can be verified using \cite[p. 223]{Kloeden1991} and \blue{\cite[Sec. 6.12]{Telatovich}}.
We note that the right hand side of the formula does not depend on $k$, and so the sampling statistics are the
same at each step.}
Each of the ODEs corresponding to these operators has explicit solutions. Using the BCH formula, and a comparison with the It\^o-Taylor expansion, we found that,
\begin{thm}
 The non-symmetric splitting scheme
 \[ \exp(D_{\dt}^{\rm{\rm II}}) \approx \exp(A) \exp(B) \]
 yields an approximation with strong order 1. The symmetric splitting scheme,
  \[ \exp(D_{\dt}^{\rm{\rm II}}) \approx \exp(\frac{A}2) \exp(B) \exp(\frac{A}2)\] gives strong order 2,
  \green{and the other symmetric splitting scheme,
  \[ \exp(\frac{B}2) \exp(A) \exp(\frac{B}2) \] has the same order \cite{Hairer2005}.}
\end{thm}
\blue{
\begin{proof}
    For the symmetric-splitting, the proof is identical to the proof of Theorem \cref{thm: order1},
    using the triangle inequality, the fact that $\exp D_{\dt}^{\rm II}$ generates a strong order 2 solution,
    and the well-known fact that the symmetric splitting for ODEs has order 2 (see \cite{Hairer2005}).

    For the non-symmetric splitting, a direct Taylor expansion shows that the low order convergence
    of the splitting reduces the overall convergence to 1. The expansion is straightforward and left
    to the interested reader.
\end{proof}
}
The non-symmetric splitting, \cref{alg:non-sym}, involves two steps, where $c_0=e^{-\Gamma\dt}$ and $c_1=\Gamma^{-1}(I-c_0)$.\green{
\begin{algorithm}
\caption{Non-symmetric splitting of truncation \rm II}
\label{alg:non-sym}
\begin{algorithmic}
\STATE{Given $x=x_0, v=v_0$, evaluate $f(x)$}
\FOR{$k=0,1,\dots,n_T-1$}
\STATE{Sample $\dw$ and $\du$ according to \eqref{eq: covtrunc2}}
\STATE{Update $v\longleftarrow c_0 v + \frac{c_1}{\dt}\left( f(x)\dt+\sigma\dw_k+\Gamma\sigma\du  \right)$}
\STATE{Update $x\longleftarrow x+\dt v-\sigma\du$}
\STATE{Evaluate $f(x)$}
\ENDFOR
\RETURN $x=x(T),v=v(T)$
\end{algorithmic}
\end{algorithm}
}
The symmetric splitting, \cref{alg:sym}, involves three steps.\green{
\begin{algorithm}
\caption{Symmetric splitting of truncation \rm II}
\label{alg:sym}
\begin{algorithmic}
\STATE{Given $x=x_0, v=v_0$, evaluate $f(x)$}
\FOR{$k=0,1,\dots,n_T-1$}
\STATE{Sample $\dw$ and $\du$ according to \eqref{eq: covtrunc2}}
\STATE{Update $x\longleftarrow x+\frac{\dt}2\left(v-\frac{\sigma \du}{\dt} \right)$}
\STATE{Update $v \longleftarrow c_0 v + \frac{c_1}{\dt}\left( f(x)\dt+\sigma\dw+\Gamma\sigma\du  \right)$}
\STATE{Update $x \longleftarrow x+\frac{\dt}2\left(v-\frac{\sigma \du}{\dt} \right)$}
\STATE{Evaluate $f(x)$}
\ENDFOR
\RETURN $x=x(T),v=v(T)$
\end{algorithmic}
\end{algorithm}
}
The numerical tests in \cref{fig: pend} and \cref{fig: LJ} confirmed the convergence orders
for truncation II with the naive and symmetric splittings methods can be found in \cref{sec: tests}.
The symmetric splitting method has the same strong convergence order as the SVV method. It is slightly easier to
implement.

\blue{
\subsection{Under-damped ($\Gamma\rightarrow 0$) and over-damped ($\Gamma \green{\gg} 1$) cases}

We consider the 2nd order symmetric splitting method above, given by the algorithm \cref{alg:sym}, in the under-damped situation, $\Gamma\rightarrow 0^+$. Assuming the fluctuation dissipation theorem \green{\eqref{eq: FDT}} holds, $\sigma\rightarrow 0^+$ as $\Gamma\rightarrow 0^+$. Then the constant \green{matrices} $c_0,c_1$ in the equations \cref{alg:sym} become, in the limit,
\begin{align}
    \lim_{\Gamma\rightarrow 0^+}c_0&=\lim_{\Gamma\rightarrow 0^+}e^{-\Gamma\dt}=I \\
    \lim_{\Gamma\rightarrow 0^+}c_1&= \lim_{\Gamma\rightarrow 0^+}\Gamma^{-1}(I-c_0)
        =\lim_{\Gamma\rightarrow 0^+}\dt e^{-\Gamma\dt}=\dt \green{I},
\end{align}
where the second limit-equality holds due to l'H\^opital's rule. As a result, the equations in \cref{alg:sym} become,
\begin{align}
    x_{n+1/2}&= x_n + \frac{\dt v_n}{2}\\
    v_{n+1}&= v_n + \dt f(x_{n+1/2})\\
    x_{n+1}&= x_{n+1/2} + \frac{\dt v_{n+1}}{2},
\end{align}
that is, the under-damped limit of the 2nd order operator splitting method is the velocity Verlet method for the corresponding deterministic Hamiltonian dynamics \green{\cite{leimkuhler2015molecular}},
\begin{equation}
    \ddot{x}=f(x).
\end{equation}

Next, we turn to the over-damped limit, $\Gamma\gg 1$ large, for truncation $\rm II$ with the non-symmetric splitting. 
 If we take $\Gamma\gg 1$ large, the term $c_0\approx 0$ is negligible and $c_1\approx\Gamma^{-1}$, so that the steps become independent of $v$:
\begin{equation}
\begin{aligned}
    v_{k+1}&=\Gamma^{-1}f(x_k)+\frac{\Gamma^{-1}\sigma\dw_k}{\dt}+\frac{\sigma\du_k}{\dt}\\
    x_{k+1}&=x_k+\dt\left( \Gamma^{-1}f(x_k)+\frac{\Gamma^{-1}\sigma\dw_k}{\dt}+\frac{\sigma\du_k}{\dt} \right)-\sigma\du_k\\
        &=x_k + \dt\Gamma^{-1}f(x_k)+\Gamma^{-1}\sigma\dw_k.
\end{aligned}
\end{equation}
In particular, the second step is the Euler-Maruyama method for the position-\green{only} first degree SDE,
\begin{equation}\label{eq: reducedSDE}
    dx = \Gamma^{-1}f(x)dt+\Gamma^{-1}\sigma dW_t.
\end{equation}
Thus the over-damped limit of the non-symmetric splitting with truncation $\rm II$ is the Euler-Maruyama method for the first degree SDE above.

Turning to the symmetric splitting in the over-damped limit, we consider the splitting
\begin{equation}
    \exp(B/2)\exp(A)\exp(B/2).
\end{equation}
In this setting, the first two steps give us the Euler-Maruyama method for the reduced first order SDE \eqref{eq: reducedSDE}. Indeed, we have
\begin{equation}
\begin{aligned}
    v_{k+1/2}&= \Gamma^{-1}f(x_k)+\frac{\Gamma^{-1}\sigma\dw_k}{\dt}+\frac{\sigma\du_k}{\dt}\\
    x_{k+1}&= x_k+ v_{k+1/2}\dt-\sigma\du_k\\
        &= x_k + \dt\left( \Gamma^{-1}f(x_k)+\frac{\Gamma^{-1}\sigma\dw_k}{\dt}+\frac{\sigma\du_k}{\dt} \right)-\sigma\du_k\\
        &= x_k + \Gamma^{-1}f(x_k)\dt+\Gamma^{-1}\sigma \dw_k,
\end{aligned}
\end{equation}
and the last equation, $x_{k+1}=x_k+\Gamma^{-1}f(x_k)\dt+\Gamma^{-1}\sigma\dw_k,$
which is independent of $v$, is precisely the Euler-Maruyama method for $dx_t=\Gamma^{-1}f(x_t)dt+\Gamma^{-1}\sigma dW_t.$
}


\subsection{Third-order truncation}\label{sec: trunc3}

Finally, we turn to the next truncation,
\begin{equation}\label{eq: trunc3}
    D_{\dt}^{\rm III}= \dt+\sum_{j=1}^n\dw^j \mathcal{X}_j + \sum_{j=1}^n\du^j [\mathcal{X}_0,\mathcal{X}_j]+ \sum_{j=1}^n\dV^j\left[[\mathcal{X}_0,\mathcal{X}_j],\mathcal{X}_0\right]
\end{equation}
where
\begin{equation}\begin{aligned}
    \dV^j&= \left[[\dt,\dw^j],\dt\right]=\frac{1}{18}\left(2I_{(0,j,0)}-2I_{(j,0,0)}+\dt I_{(j,0)}-\dt I_{(0,j)}\right)\\
        &=\frac{1}{9}\left( I_{(0,j,0)}-I_{(j,0,0)}-\dt\du^j \right)
\end{aligned}\end{equation}
and $\dV=(\dV^1,\dots,\dV^n)$.
This comes from \cite[eq. 2.15]{Misawa2000}, using the fact that $[[\mathcal{X}_0,\mathcal{X}_i],\mathcal{X}_j]=0$ when $i,j\neq 0$ for additive noise.
Since we have the following expressions for all relevant brackets,
\begin{equation}\begin{aligned}
    \mathcal{X}_0&= v\cdot\nabla_x+(f(x)-\Gamma v)\cdot\nabla_v\\
    \mathcal{X}_j&= \sigma^j\cdot\nabla_v\\
    \left[\mathcal{X}_0,\mathcal{X}_j\right]&= -\sigma^j\cdot \nabla_x+\Gamma\sigma^j\cdot \nabla_v\\
    \left[[\mathcal{X}_0,\mathcal{X}_j],\mathcal{X}_0\right]&= \Gamma\sigma^j\cdot \nabla_x-( Df(x)\sigma^j+\Gamma^2\sigma^j)\cdot\nabla_v,
\end{aligned}\end{equation}
we can rewrite $D_{\dt}^{\rm III}$ in \eqref{eq: trunc3} as
\begin{equation}\begin{aligned}
    D_{\dt}^{\rm III}&= \left(v\dt-\sigma\du+\Gamma\sigma\dV\right)\cdot\nabla_x\\
                    &+ \left( (f(x)-\Gamma v)\dt+\sigma\dw+\Gamma\sigma\du-(Df(x)\sigma+\Gamma^2\sigma)\dV \right)\cdot\nabla_v,
\end{aligned}\end{equation}
which gives us the ODEs
\begin{equation}\begin{aligned}\begin{cases}
 x'&= v\dt-\sigma\du+\Gamma\sigma\dV\\
 v'&=(f(x)-\Gamma v)\dt+\sigma\dw+\Gamma\sigma\du-(Df(x)\sigma+\Gamma^2\sigma)\dV.
 \end{cases}\label{eq: III}
\end{aligned}\end{equation}
\blue{We do not have a complete proof that $\exp D_{\dt}^{\rm III}$ generates solutions with strong order 3. Although
we have the intuition that  this can be proven in a similar way to the way we proved the convergence
orders of the solutions generated by $\exp D_{\dt}^{\rm I}$ and $\exp D_{\dt}^{\rm II}$, but the calculations
will be more complicated. The numerical evidence supports this conjecture, see \cref{fig: pend} and \cref{fig: LJ}.}

For the numerical implementation, we use the splitting $D_{\dt}^{\rm III}=A+B$ where,
\begin{equation}\label{eq: AB-III}
 \begin{aligned}
    A&= \left( v\dt-\sigma \du+\Gamma\sigma \dV \right)\cdot\nabla_x\\
    B&= \left((f(x)-\Gamma v)\dt+\sigma\dw +\Gamma\sigma\du-(Df(x)\sigma+\Gamma^2\sigma)\dV\right)\cdot\nabla_v,\\
    \end{aligned}
\end{equation}
and the two corresponding ODEs are given by
\begin{equation}\begin{aligned}\begin{cases}
 x'&= v\dt-\sigma\du+\Gamma\sigma\dV\\
 v'&=0
\end{cases}\end{aligned}\end{equation}
and
\begin{equation}\begin{aligned}\begin{cases}
 x'&= 0\\
 v'&=(f(x)-\Gamma v)\dt+\sigma\dw+\Gamma\sigma\du-(Df(x)\sigma+\Gamma^2\sigma)\dV.
\end{cases}\end{aligned}\end{equation}
They have explicit solutions\blue{, at time $t_{k+1}=(k+1)\dt$, given by,}
\begin{equation}\begin{aligned}
x_{k+1}&=x_k+\dt\left(v_k-\frac{\sigma \du_k}{\dt}+\frac{\Gamma\sigma\dV_k}{\dt}\right)
\end{aligned}\end{equation}
and
\begin{equation}\begin{aligned}
v_{k+1}&=c_0 v_k + \frac{c_1}{\dt}\left( f(x_k)\dt+\sigma\dw_k+\Gamma\sigma\du_k-(\nabla_x f(x_k)\sigma+\Gamma^2\sigma)\dV_k \right)
\end{aligned}\end{equation}
where again $c_0=\exp(-\Gamma \dt)$ and $c_1=\Gamma^{-1}(I-c_0)$.
If \green{the coefficient $\nabla_x f(x_k)\sigma + \Gamma^2\sigma$ of $\dV_k$ is zero}, this method is reduced to truncation II. \blue{We acknowledge and highlight that truncation $\rm III$ requires computing the derivative of the potential function $f$. However, this is one fewer derivative than is required by the It\^o-Taylor approximation, which requires the second order derivatives. (See the last formula in \eqref{eq: IT0}.)}

In the implementation, the term $\nabla_x f(x)\sigma\green{\dV}$ will be approximated by a finite-difference formula,
\begin{equation}
  \nabla_x f(x)\sigma \Delta V \approx \frac{f(x+\varepsilon \sigma \Delta V) - f(x)}{\varepsilon},\text{\quad}0<\epsilon\green{\ll} 1.
\end{equation}
\blue{In our numerical convergence tests, we took $\epsilon$ to be equal to the small step size $\delta{t}=2^{-19}$ for the discretized Brownian motion, which stays fixed while the larger step size $\Delta{t}$ for the discretized solution of the SDE varies.}

In principle, the symmetric splitting methods applied to ODEs have order 2. In order to achieve higher order of accuracy, we solve the ODEs \eqref{eq: III} using the Neri's splitting method \cite{Yoshida1990}, which consists of alternating the operators in \eqref{eq: AB-III} three times:
\begin{equation}\begin{aligned}
    \exp(D^{\rm III})&\approx \exp(a_1A)\exp(b_1B)\exp(a_2A)\exp(b_2B)\exp(a_3A)\exp(b_3B)\exp(a_4A),\\
    \text{where\quad}
        a_1&=\frac{1}{2(2-\sqrt[3]{2})},\text{\quad}a_2=\frac{1}{2}-a_1,\text{\quad}a_3=a_2,\text{\quad}a_4=a_1,\\
    \text{and\quad}
        b_1&=\frac{1}{2-\sqrt[3]{2}},\text{\quad}b_2=1-2b_1,\text{\quad}b_3=b_1.
\end{aligned}\end{equation}
\blue{We wish to highlight that to obtain such high order of accuracy at least one of the steps has to be negative, which is well known
(see \cite[Ch. 3]{Hairer2005}).}



\bigskip

For practical implementations,  the joint covariances of $\dw_k, \du_k$, and $\dV_k$ are needed in order to sample these mean-zero Gaussian random variables.
Using \cite[p. 223]{Kloeden1991}, it can be shown (see \blue{\cite[Sec. 6.12]{Telatovich}}) that
\begin{equation}\begin{aligned}\label{eq: covtrunc3}
C\overset{\text{def}}{=}     &\left( \begin{matrix} \mathbb{E}(\dw_k^2)&\mathbb{E}(\dw_k\du_k)&\mathbb{E}(\dw_k\dV_k)\\
                            \mathbb{E}(\du_k\dw_k)&\mathbb{E}(\du_k^2)&\mathbb{E}(\du_k\dV_k)\\
                            \mathbb{E}(\dV_k\dw_k)&\mathbb{E}(\dV_k\du_k)&\mathbb{E}(\dV_k^2)\end{matrix} \right)\\
    &= \left(\begin{matrix} \dt I_{n\times n}& 0_{n\times n}& 0_{n\times n}\\ 
                            0_{n\times n}&\frac{\dt^3}{12}I_{n\times n} & \frac{-\dt^4}{216}I_{n\times n                       }\\
                            0_{n\times n}& \frac{-\dt^4}{216}I_{n\times n} & \frac{\dt^5}{2430}\dt^5 I_{n\times n} \end{matrix}\right).
\end{aligned}\end{equation}
To sample $(\dw_k,\du_k,\dV_k)$, we use the Cholesky decomposition $C=LL^T,$ and sample the noise by multiplying $L^T$ to independent Gaussian random variables.


\section{\green{Proofs of the convergence theorems}}\label{sec: proofs}

\subsection{The analysis of the stochastic velocity Verlet method}\label{sec: svv}

Here we give a detailed proof of the following:

\begin{thm}
  The SVV algorithm has strong order 2.
\end{thm}

\begin{proof}
We start with the displacement component. We compare
\begin{equation}
    x_{k+1}=x_k+v_k\dt+(f(x_k)-\Gamma v_k)\frac{\dt^2}{2}+\sum_{j=1}^n\sigma^j I_{(j,0),k}
\end{equation}
i.e., the It\^o-Taylor approximation with strong order 2, with the SVV method,
\begin{equation}\begin{aligned}
      \tilde{x}_{k+1}  &=x_k+v_k\dt+(f(x_k)-\Gamma v_k)\frac{\dt^2}{2}+\sum_{j=1}^n\sigma^j I_{(j,0),k}\\
            &+\sum_{j=1}^n \int_{t_k}^{t_{k+1}} \int_{t_k}^t \left(e^{-\Gamma(t-s)}-I\right) \sigma^jdW_s^j dt
                +\mathcal{O}(\dt^3).
\end{aligned}\end{equation}
The remainder term $R_k, k=0,\dots,n_T-1$ after the $k$th step is a discrete martingale \cite[pg. 195]{Kloeden1991}
\begin{equation}\begin{aligned}
R_k&= \sum_{j=1}^n \int_{t_k}^{t_{k+1}} \int_{t_k}^t \left(e^{-\Gamma(t-s)}-I\right)\sigma^j dW_s^j dt,
\end{aligned}\end{equation}
and it remains to show that
\begin{equation}\label{eq: conv}
    \mathbb{E}\left( \max_{1\leq m\leq n_T}\left| \sum_{k=0}^{m-1} R_k \right|^2\right)\leq C\dt^4.
\end{equation}
 To that end, first notice that $\mathbb{E}(R_k^T R_l)=\delta_{kl}\mathbb{E}(|R_k|^2)$ (where $\delta_{kl}=1$ if $k=l$ and 0 otherwise), since the increments of the Wiener process are independent. Next, since $R_k$ is a martingale, we may apply the discrete version of Doob's lemma with $p=2$\cite[eq. 2.3.7]{Kloeden1991} to obtain an estimate for \eqref{eq: conv}:
\begin{equation}\begin{aligned}
\mathbb{E}\left(\max_{1\leq m\leq n_T}\left| \sum_{k=0}^{m-1} R_k \right|^2\right)
    &\leq 4 \mathbb{E}\left( \left| \sum_{k=0}^{n_T-1} R_k \right|^2 \right)\text{\quad by Doob's lemma}\\
    &\leq 4\sum_{k=0}^{n_T-1} \mathbb{E}(R_k^T R_k)\text{\quad since }\mathbb{E}(R_k^T R_l)=\delta_{kl}\mathbb{E}(\blue{R_k^T R_k})\\
    &= 4n_T \mathbb{E}(R_k^T R_k).
\end{aligned}\end{equation}
Notice $e^{-\Gamma(t-s)}-I_n=\mathcal{O}(\dt)$ for $t_k\leq s\leq t\leq t_{k+1}$, and recall that
$\mathbb{E}(I_{(j,0)}^2)=\mathcal{O}(\dt^3)$ for all $j=1,2,\cdots,n$ \cite[pg. 172, exercise 5.2.7]{Kloeden1991}.
As $\mathbb{E}(I_{(j_1,0)}I_{(j_2,0)})=0$ for distinct $j_1,j_2=1,2,\cdots,n$ the random variables
$W^{j_1},W^{j_2}$ are independent, see \cite[p. 223, eq. 5.12.7]{Kloeden1991} for more details), we have
\begin{equation}\begin{aligned}
\mathbb{E}(R_k^T R_k)&=\mathbb{E}\left((  \mathcal{O}(\dt^2) \left(\sum_{j=1}^n \sigma^{j}I_{(j,0),k}\right)^T \left(\sum_{l=1}^n \sigma^{j}I_{(j,0),l}\right)  \right)\\
        &= \mathcal{O}(\dt^2)\mathbb{E}\left(\sum_{j=1}^n I_{(j,0),k}^2 \left(\sigma^j\right)^T\sigma^j \right)\\
        &= \mathcal{O}(\dt^2)\sum_{j=1}^n \mathbb{E}\left( I_{(j,0),k}^2 \right)\left( \sigma^j \right)^T \sigma^j \\
        &= \mathcal{O}(\dt^2)\mathcal{O}(\dt^3)\text{\quad since }n<<n_T\\
        &=\mathcal{O}(\dt^5).
\end{aligned}\end{equation}
Therefore, $\mathbb{E}\left( \max_{1\leq m\leq n_T}\left| \sum_{k=0}^{m-1} R_k \right|^2 \right)\leq 4n_T \mathcal{O}(\dt^5)=\mathcal{O}(\dt^4)=\mathcal{O}(\dt^{2\gamma})$, so convergence criterion  is satisfied for $\gamma=2$.

For the velocity components, with $c_0=e^{-\Gamma\dt}$, $c_1=(\Gamma\dt)^{-1}(I-c_0)$, and $c_2=\frac{1}{\dt^2}\int_0^\dt \Gamma^{-1}(I-e^{-\Gamma t})dt$, where $I$ is the $n\times n$ identity matrix, we compare the order 2 It\^{o} Taylor approximation (recall \eqref{eq: IT0}),
\begin{equation}\begin{aligned}
    v_{k+1}&=v_k+(f(x_k)-\Gamma v_k)\dt+(\nabla_x f(x_k)v_k-\Gamma(f(x_k)-\Gamma v_k))\frac{\dt^2}{2}\\
        &+\sum_{j=1}^n\sigma^j\dw_k^j -\sum_{j=1}^n\Gamma \sigma^j I_{(j,0),k}
\end{aligned}\end{equation}
with the SVV algorithm,
\begin{equation}\begin{aligned}
 \tilde{v}_{k+1} =c_0v_k&+c_1 f(x_k)\dt +c_2 (f(x_{k+1})-f(x_k))\dt+ \sum_{j=1}^n \int_{t_k}^{t_{k+1}}  e^{-\Gamma(\dt-s)}\sigma^j dW_s^j\\
    = v_k&+(f(x_k)-\Gamma v_k)\dt +[\nabla_x f(x_k)v_k-\Gamma(f(x_k)-\Gamma v_k)]\frac{\dt^2}{2}+\sum_{j=1}^n\sigma^j\dw^j\\
    &-\sum_{j=1}^n\Gamma\sigma^j I_{(j,0)}+\sum_{j=1}^n \int_{t_k}^{t_{k+1}}  \left(e^{-\Gamma(\dt-s)}-I+\Gamma(\dt-s)\right)\sigma^j dW_s^j,
\end{aligned}\end{equation}
where the last equality is non-trivial and holds from Taylor expanding the $c$'s and $f(x_{k+1})-f(x_k)$.
We see that all the terms match up to the remainder terms
\begin{equation}\begin{aligned}
R_k=\sum_{j=1}^n \int_{t_k}^{t_{k+1}}  \left(e^{-\Gamma(\dt-s)}-I+\Gamma(\dt-s)\right)\sigma^j dW_s^j,
    \text{\quad} k=0,\dots,n_T-1,
\end{aligned}\end{equation}
and notice that $R_k$ is a discrete martingale. Once again, $\mathbb{E}\left(R_k^T R_l\right)=\delta_{kl}\mathbb{E}\left(R_k^T R_k\right)$ since Brownian increments are independent and the It\^{o} integrals are taken over disjoint intervals. Then we see that
\begin{equation}\begin{aligned}
\mathbb{E}\left( R_k^T R_k \right)&= \int_{t_k}^{t_{k+1}}  (e^{-\Gamma(\dt-s)}-I+\Gamma(\dt-s))^2 ds
    \text{\quad by It\^{o}'s isometry \cite[Cor. 3.7]{Oksendal1995}}\\
        &= C_k \int_0^\dt s^4 ds \text{\quad by Taylor series approximation}\\
        & = \frac{C_k}{5}\dt^5
\end{aligned}\end{equation}
for some constant $C_k>0$, which implies that,
\begin{equation}\begin{aligned}
\mathbb{E}\left( \max_{1\leq m\leq n_T}\left|\sum_{k=0}^{m-1}R_k\right|^2 \right)
    &\leq 4 \mathbb{E}\left( \left|\sum_{k=0}^{n_T-1} R_k^T R_k \right| \right)\text{\quad by discrete Doob's lemma}\\
    &\leq 4 \mathbb{E}\left( \sum_{k=0}^{n_T-1} R_k^T R_k \right) \text{\quad since }\mathbb{E}(R_k^T R_l)=\delta_{kl}\mathbb{E}\left(\blue{R_k^T R_k}\right)\\
    &\leq\frac{4C}{5}\sum_{k=0}^{n_T-1} \mathbb{E}\left(R_k^T R_k\right)\text{\quad where }C:=\max_{0\leq k<n_T}C_k\\
    &=\frac{4C}{5}n_T \mathcal{O}(\dt^5)\\
    &= \mathcal{O}(\dt^4)\text{\quad since }n_T\dt= T.
\end{aligned}\end{equation}
Therefore convergence criterion \eqref{eq: conv} is satisfied.
\end{proof}



\subsection{Order 1 convergence of the direct splitting method}\label{sec: pfdirectsplit}

\blue{
\begin{thm}
    The direct splitting method produces a strong order 1 scheme.
\end{thm}
\begin{proof}
We recall that the direct splitting is,
\begin{equation}\begin{cases}
    x' &= v\\
    v' &= 0
\end{cases}\text{\quad}
\begin{cases}
    x' &= 0\\
    v' &= f(x)-\Gamma v +\sigma W'.
\end{cases}
\end{equation}
Solving the first system for one step $\dt$ with $v$ constant, we get,
\begin{equation}
    x(\dt)=x+v\dt = x + g_{(0)}(x)\dt,
\end{equation}
so that the $x$ component agrees with the order 1 It\^o-Taylor expansion.
The second system, where $x$ is constant, can be viewed as a first order linear system of SDEs,
\begin{equation}
    dv_t = \left( -\Gamma v_t + f(x) \right)dt + \sigma dW_t,
\end{equation}
or more abstractly,
\begin{equation}
    dv_t = \left( Av_t + B \right)dt + CdW_t,
\end{equation}
where $A$ and $C$ are constant $d\times d$ matrices, $B\in\R^d$ is a constant vector, and $W_t$ is $d$-dimensional
Wiener process. The notation in the abstract form is more convenient for solving the equations.
It is well known how to solve such an SDE (see \cite[Example 5.3]{Oksendal1995} for a similar SDE)
but we present the worked solution for completeness.
We multiply $v_t$ by an integrating factor $\exp(-At)$ and use the It\^o formula:
\begin{equation}\begin{aligned}
    d\left( \exp(-At)v_t \right)
        &= -A\exp(-At)v_t dt + \exp(-At)dv_t \\
        &= -A\exp(-At)v_t dt + \exp(-At)\left( Av_t dt + B dt + C dW_t \right)\\
        &= \exp(-At)B dt + \exp(-At)C dW_t\text{\quad since }A\exp(-At)=\exp(-At)A,
\end{aligned}\end{equation}
which is just a formal way of writing the integral equation,
\begin{equation}
    \exp(-At)v_t = v_0 + \int_0^t \exp(-As)B ds + \int_0^t \exp(-As)CdW_s,
\end{equation}
that is,
\begin{equation}
    v_t =\exp(At)\left[ v_0 + A^{-1}(I-\exp(-At))B+\int_0^t \exp(-As)CdW_s \right].
\end{equation}
Now apply the formula above to the original equation, where $A=-\Gamma$, $B=f(x)$, and
$C=\sigma$:
\begin{equation}
    v_t = \exp(-\Gamma t)\left[ v_0 - \Gamma^{-1}(I-\exp(-\Gamma t))f(x)+\int_0^t \exp(\Gamma s)\sigma dW_s \right],
\end{equation}
that is,
\begin{equation}
    v_t = \exp(-\Gamma t)v_0 +\Gamma^{-1}(I-\exp(-\Gamma t))f(x)+\int_0^t \exp(-\Gamma(t-s))\sigma dW_s,
\end{equation}
or $v_t = c_0(t) v_0 + c_1(t) f(x) + c_2(t)$, where
\begin{equation}\begin{aligned}
    c_0(t) &= \exp(-\Gamma t) \\
    c_1(t) &= \Gamma^{-1}(I-\exp(-\Gamma t))f(x) \\
    c_2(t) &= \int_0^t \exp(-\Gamma(t-s))\sigma dW_s.
\end{aligned}\end{equation}
Setting $v_0=v$ and $t=\dt$, and assuming $\dt$ is much less than the spectral radius of $\Gamma^{-1}$, we observe:
\begin{equation}\begin{aligned}
    c_0 v &= v - \Gamma v \dt + \mathcal{O}(\dt^2) \text{\quad linear approximation}\\
    c_1 f(x) &= f(x)\dt + \mathcal{O}(\dt^2)\text{\quad linear approximation} \\
    c_2 &= \sigma\dw - \int_0^{\dt}\Gamma \exp(-\Gamma(\dt-s))\sigma W_s ds,
\end{aligned}\end{equation}
where the last (exact) formula holds by the stochastic integration by parts formula (see \cite[Theorem. 4.5]{Oksendal1995}).
The stochastic integral has $\mathcal{O}(\dt^{3/2})$, since its variance is $\mathcal{O}(\dt^3)$
(see \cite{Kloeden1991}). A constant approximation of the integrand
$\Gamma\exp(-\Gamma(\dt-s))\sigma = \Gamma\sigma+\mathcal{O}(s)$ gives us $g_{(j,0)}(v)$ from the It\^o-Taylor
expansion:
\begin{equation}\begin{aligned}
    -\int_0^{\dt}\Gamma\exp(-\Gamma(\dt-s))\sigma W_s ds &= -\Gamma\sigma \int_0^{\dt}W_s ds\green{+\dt\mathcal{O}(\dt^{3/2})}\\
    &= g_{(j,0)}(v)I_{(j,0)}+\mathcal{O}(\dt^{5/2}).
\end{aligned}\end{equation}
Therefore, using linear approximation, we obtain,
\begin{equation}\begin{aligned}
    v(\dt) &= c_0 v + c_1 f(x) + c_2 \\
        &= v - \Gamma v\dt + f(x)\dt + \sigma\dw - \Gamma\sigma I_{(J,0)}+\mathcal{O}(\dt^2) \\
        &= v + g_{(0)}(v)\dt + g_{(J)}(v)I_{(J)} + \sum_{j=1}^n g_{(j,0)}(v)I_{(j,0)}+ \mathcal{O}(\dt^2).
\end{aligned}\end{equation}
That is, the $v$ component agrees with the order 1 It\^o-Taylor expansion. Therefore the direct splitting method
produces strong order 1 schemes. They fail to attain higher order convergence because of the absence of
higher order stochastic integrals $I_{(j,0)}$ which are present in the order 2 It\^o-Taylor expansion.

The fact that we get $g_{(j,0)}$ terms gives us reason to believe that the direct splitting
can attain strong order 2 convergence in the $v$ component. Indeed, it may be that a quadratic approximation
of $c_0v$ and $c_1 f(x)$ would give us $g_{(0,0)}(v)$. However, the low order convergence of the
$x$ component restricts us from obtaining higher order convergence over
all coordinates.
\end{proof}
 }

\subsection{First-order convergence of $\exp D_\dt^{\rm I}(x,v)$}\label{sec: pftrunc1}

Abbreviating $z=(x,v)$, we prove that $\exp D_{\dt}^{\rm I}(z)$ matches the strong order $\gamma=1$ It\^o-Taylor
expansion up to terms of order $\dt^{1.5}$. We borrow much of the notation from \cref{sec: IT}.
\begin{proof}
To simplify calculations, we will use the obvious facts that
\begin{equation}
    \mathcal{X}_j(\phi)=L^j(\phi)
\end{equation}
for any $C^1$ function $\phi$, and
\begin{equation}
    \mathcal{X}_0(\phi)=L^0(\phi)
\end{equation}
for any function $\phi$ that is linear in $v$ (for such functions $\partial^2/\partial v^2$ vanishes identically).
Recall that
\begin{equation}
    D_{\dt}^{\rm I}= \dt\X_0+\sum_{j=1}^m \dw^j \X_j.
\end{equation}
Since $z=(x,v)$ is linear in $v$,
\begin{equation}
    D_{\dt}^{\rm I}(z)=\dt L^0(z)+\sum_{j=1}^m \dw^j L^j (z)=g_{(0)}I_{(0)}+\sum_{j=1}^m g_{(j)}I_{(j)}
\end{equation}
where the $g_\alpha$ and $I_\alpha$ are the coefficient functions and stochastic integrals, respectively, in the It\^o-Taylor expansions \cref{eq: IT0}(we often suppress the arguments $g=g(z)$).
Then, since $g_{(0)}(z)=(v, f(x)-\Gamma v)^T$ and $g_{(j)}(z)=(0,\sigma^j)^T$ are linear in $v$,
\begin{equation}
    \begin{aligned}
        \frac{1}{2}\left( D_{\dt}^{\rm I} \right)^2(z)
            &= \frac{1}{2}D_{\dt}^{\rm I}\left( g_{(0)}(z)\dt+\sum_{j=1}^m g_{(j)}(z)I_{(j)} \right)\\
            &= \frac{\dt}{2} D_{\dt}^{\rm I}g_{(0)}(z)+\sum_{j=1}^m \frac{\dw^j}{2}D_{\dt}^{\rm I} g_{(j)}(z) \\
            &= \frac{\dt}{2} D_{\dt}^{\rm I}g_{(0)}(z)\text{\quad since $g_{(j)}(z)$ is constant}\\
            &= \frac{\dt}{2}\left( \dt g_{(0,0)}(z)+\sum_{j=1}^m \dw^j g_{(j,0)}(z) \right) \\
            &= g_{(0,0)}(z)I_{(0,0)} + \sum_{j=1}^m \frac{\dt\dw^j}{2}g_{(j,0)}(z).
    \end{aligned}
\end{equation}
Combining the expressions for $D_{\dt}^{\rm I}(z)$ and $\frac{1}{2}\left(D_{\dt}^{\rm I}\right)^2(z)$ above,
we therefore have,
\begin{equation}
\begin{aligned}
    \left( \exp D_{\dt}^{\rm I} \right)(z)
        &= z + D_{\dt}^{\rm I}(z)+\frac{1}{2}\left( D_{\dt}^{\rm I} \right)^2(z)+H.O.T. \\
        &= z + g_{(0)}(z)I_{(0)} + \sum_{j=1}^m g_{(j)}(z)I_{(j)}\\
            &+g_{(0,0)}(z)I_{(0,0)} + \sum_{j=1}^m g_{(j,0)}(z) \frac{\dt\dw^j}{2}+H.O.T.,
\end{aligned}
\end{equation}
where the higher order terms (denoted by H.O.T.) do not involve $I_{(j,0)}$. This is because $I_{(j,0)}$
cannot be written as a multiple of $\dw^j$. Indeed, for any real constant $c$, $I_{(j,0)}-c\dw^j$ has non-zero variance.
Comparing with the strong order $\gamma=1,2$ It\^o-Taylor expansions,
\begin{equation}
\begin{aligned}
    z_{k+1}^{\gamma =1}&= z_k + g_{(0)}(z_k)I_{(0)} + \sum_{j=1}^m g_{(j)}(z_k)I_{(j)}\text{\quad order }\gamma=1 \\
    z_{k+1}^{\gamma = 2}&= z_{k+1}^{\gamma = 1} + g_{(0,0)}(z_k)I_{(0,0)}+\sum_{j=1}^m g_{(j,0)}(z_k)I_{(j,0)}
            \text{\quad order }\gamma=2,
\end{aligned}
\end{equation}
we see that $\exp D_{\dt}^{\rm I}(z)$
converges with strong order 1 but not 2 because $\exp D_{\dt}^{\rm I}(z)$ lacks $I_{(j,0)}$.
\end{proof}


\subsection{Second-order convergence of $\exp D_\dt^{\rm II}(x,v)$}\label{sec: pftrunc2}

We prove the second order convergence of $\exp D_{\dt}^{\rm II}(z)$ where $z=(x,v)$, by comparison to the
It\^o-Taylor expansions we derived in \cref{sec: IT} (the section explains the notation $g_\alpha, I_\alpha$).
\begin{proof}
We recall that,
\begin{equation}
    D_{\dt}^{\rm II}= \dt\X_0 + \sum_{j=1}^m \dw^j \X_j+\frac{1}{2}\sum_{j=1}^m[\dt,\dw^j][\X_0,\X_j],
\end{equation}
where $[\X_0,\X_j]=\X_0\X_j-\X_j\X_0$ and
\begin{equation}
    \frac{1}{2}[\dt,\dw^j]= \frac{\dt\dw^j}{2}-I_{(j,0)}.
\end{equation}
We know from \cref{sec: pftrunc1} that $\dt\X_0(z)=g_{(0)}(z)I_{(0)}$ and $\dw^j\X_j=g_{(j)}(z)I_{(j)}$. Now, since $g_{(j)}(z)=(0,\sigma^j)^T$
is constant, $\X_0\X_j(z)=\X_0 g_{(j)}(z)=0$, and therefore
\begin{equation}
    [\X_0,\X_j](z)=-\X_j\X_0(z)=-\X_j g_{(0)}(z)=-L^j g_{(0,)}(z)=-g_{(j,0)}(z).
\end{equation}
It follows that the linear term in $\exp D_{\dt}^{\rm II}(z)$ is,
\begin{equation}
\begin{aligned}
    D_{\dt}^{\rm II}(z) &= g_{(0)}(z)I_{(0)}+\sum_{j=1}^m g_{(j)}(z)I_{(j)}
        +\sum_{j=1}^m \left( \frac{\dt\dw^j}{2}-I_{(j,0)} \right)\left( -g_{(j,0)}(z) \right)\\
        &= g_{(0)}(z)I_{(0)}+\sum_{j=1}^m g_{(j)}(z)I_{(j)}+\sum_{j=1}^m g_{(j,0)}(z)I_{(j,0)}
            -\sum_{j=1}^m \frac{\dt\dw^j}{2}g_{(j,0)}(z).
\end{aligned}
\end{equation}
For the quadratic term $\frac{1}{2}\left( D_{\dt}^{\rm II} \right)^2(z)$, we compute (suppressing the argument $z$ of the functions $g_{(\cdot)}$),
\begin{equation}
\begin{aligned}
    \frac{1}{2}\left(D_{\dt}^{\rm II}\right)^2(z)
        &= \frac{1}{2}D_{\dt}^{\rm II}\left( g_{(0)}I_{(0)}+\sum_{j=1}^m g_{(j)}I_{(j)}+\sum_{j=1}^m g_{(j,0)}I_{(j,0)}
            -\sum_{j=1}^m \frac{\dt\dw^j}{2}g_{(j,0)} \right)\\
        &= \frac{\dt}{2}D_{\dt}^{\rm II}g_{(0)}+\sum_{j=1}^m \frac{\dw^j}{2}D_{\dt}^{\rm II}g_{(j)}\\
        &+\sum_{j=1}^m I_{(j,0)}D_{\dt}^{\rm II}g_{(j,0)}-\sum_{j=1}^m\frac{\dt\dw^j}{2}D_{\dt}^{\rm II}g_{(j,0)}\\
        &= \frac{\dt}{2}\left( \dt\X_0 g_{(0)}+\sum_{j=1}^m \dw^j \X_j g_{(0)} \right)
        \text{since $g_{(j)}$ and $g_{(j,0)}$ are constant}\\
        &= \frac{\dt^2}{2}g_{(0,0)}+\sum_{j=1}^m \frac{\dt\dw^j}{2}g_{(j,0)}.
\end{aligned}
\end{equation}
Approximating $\exp D_{\dt}^{\rm II}(z)$ by the quadratic $z+D_{\dt}^{\rm II}(z)+\frac{1}{2}\left(D_{\dt}^{\rm II}\right)^2(z)$,
we see that the terms $\dt\dw^j/2$ cancel, and we obtain,
\begin{equation}
\begin{aligned}
    \left( \exp D_{\dt}^{\rm II} \right)(z)
        &= z + g_{(0)}(z)I_{(0)}+\sum_{j=1}^m g_{(j)}(z)I_{(j)}+g_{(0,0)}(z)I_{(0,0)}\\
        &+\sum_{j=1}^m g_{(j,0)}(z)I_{(j,0)}
            +H.O.T.,
\end{aligned}
\end{equation}
which equals the strong order 2 It\^o-Taylor expansion up to the higher order terms.
\end{proof}
We believe it can be shown that $\exp D_{\dt}^{\rm III}(z)$ converges with strong order 3,
using a similar method, but we have yet to do this.


\section{Numerical tests}\label{sec: tests}

\subsection{A one-dimensional pendulum model}
\begin{figure}[htbp]
\begin{center}
\includegraphics[scale=.16]{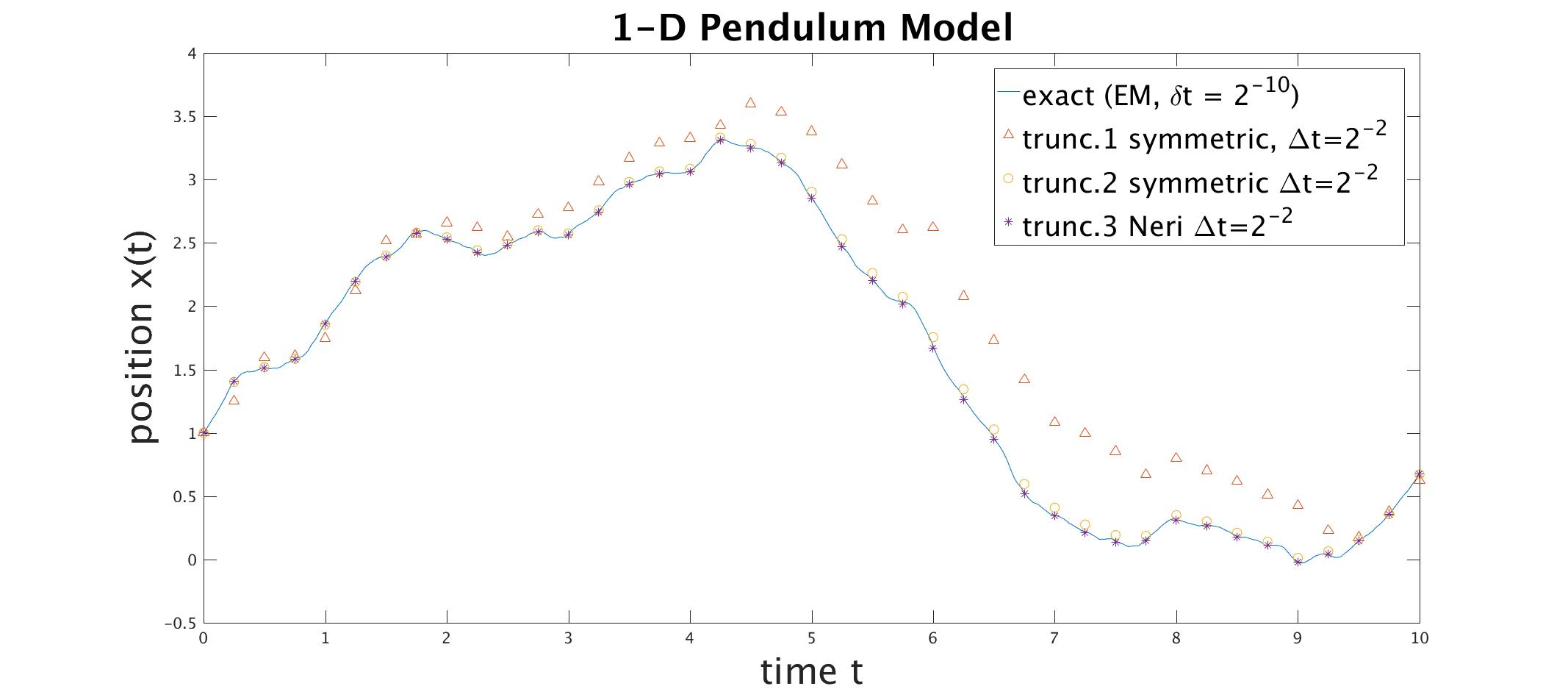}
\caption{A trajectory for the simple 1-d pendulum model, using different splitting methods\blue{, but the same reference trajectory of the driving Brownian motion.}}
\label{fig: path}
\end{center}
\end{figure}

In the first \blue{set} of experiments, we consider the one-dimensional  pendulum model:
\begin{equation}
f(x)=-\sin(x)
\end{equation}
when $x\in \R$. We start with a simulation of one trajectory using truncation I (with the symmetric splitting), truncation II (symmetric splitting) and truncation II (Neri splitting). In \cref{fig: path}, using step size $\Delta{t}=2^{-2}$, we compare them with the solution generated by the Euler Maruyama method with smaller step size $\delta{t}=2^{-10}$, \green{subsequently} viewed as the `exact' solution. All the solutions are generated from the same realization of the Brownian motion. Clearly the accuracy improves as the truncation number increases.

\begin{figure}[htbp]
\begin{center}
\includegraphics[scale=.18]{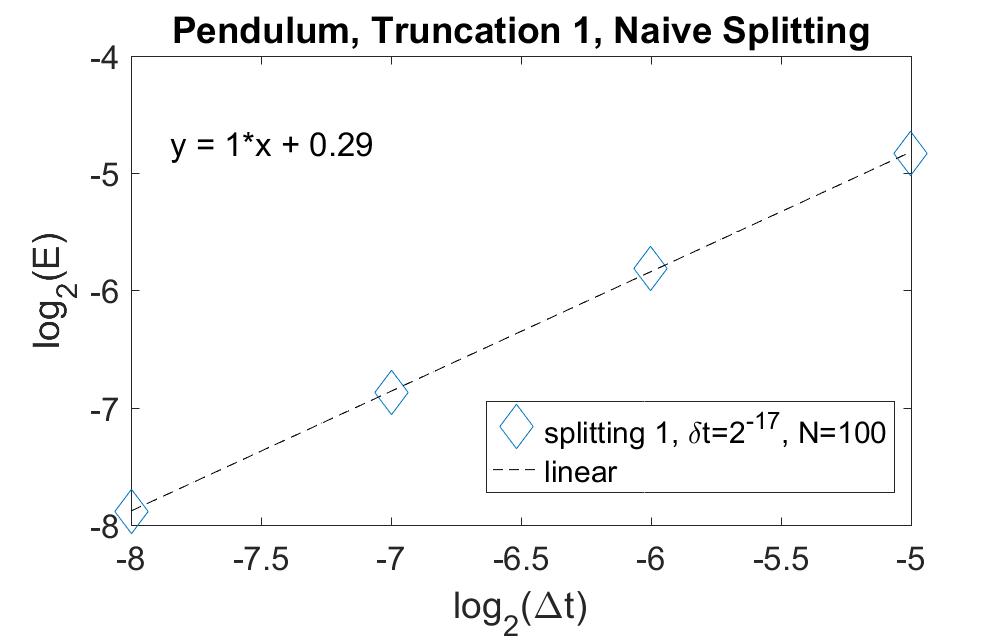}
\includegraphics[scale=.18]{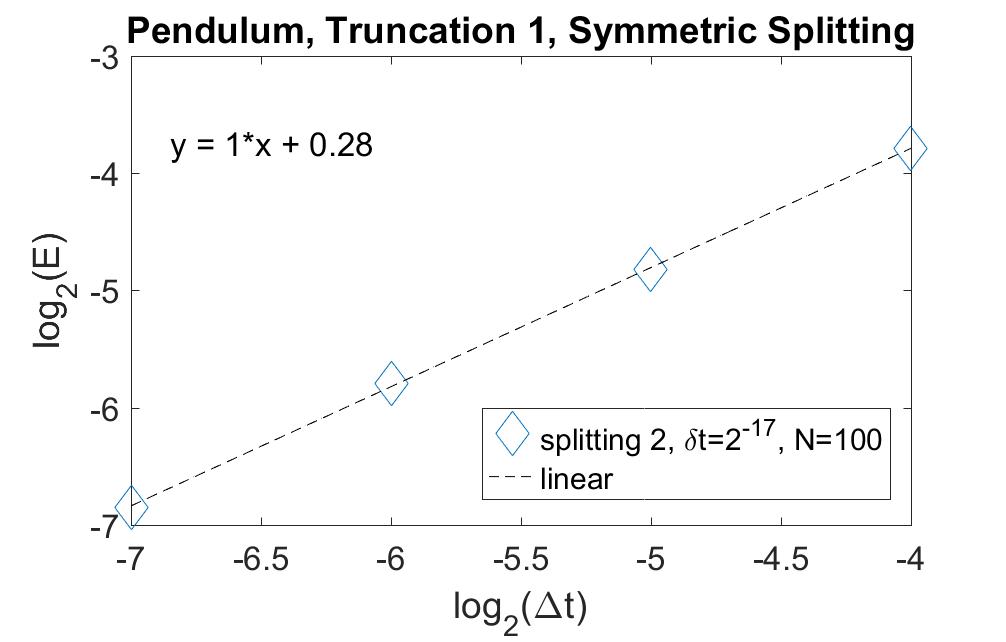}
\includegraphics[scale=.18]{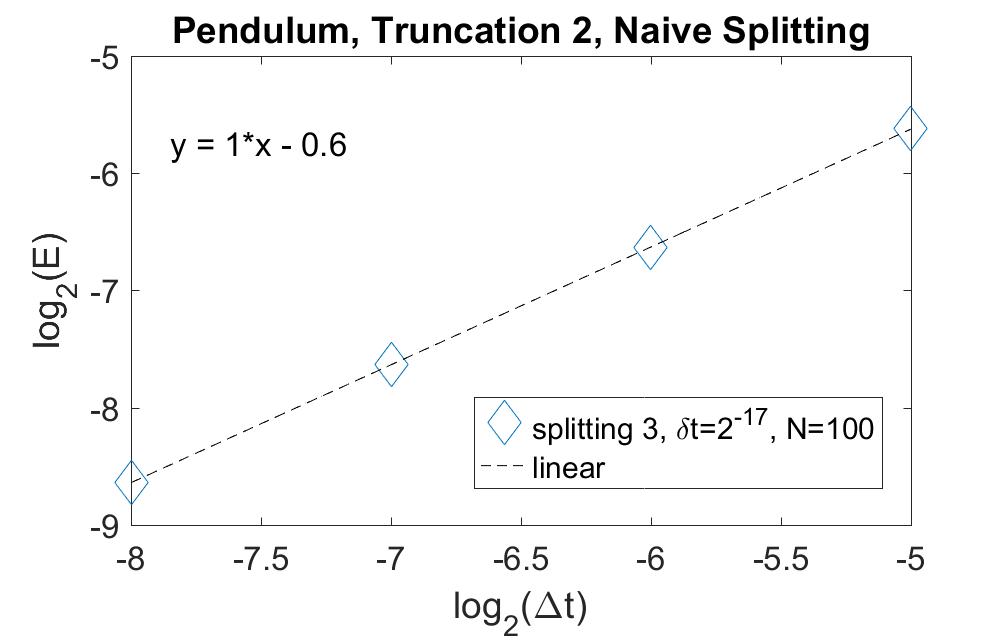}
\includegraphics[scale=.18]{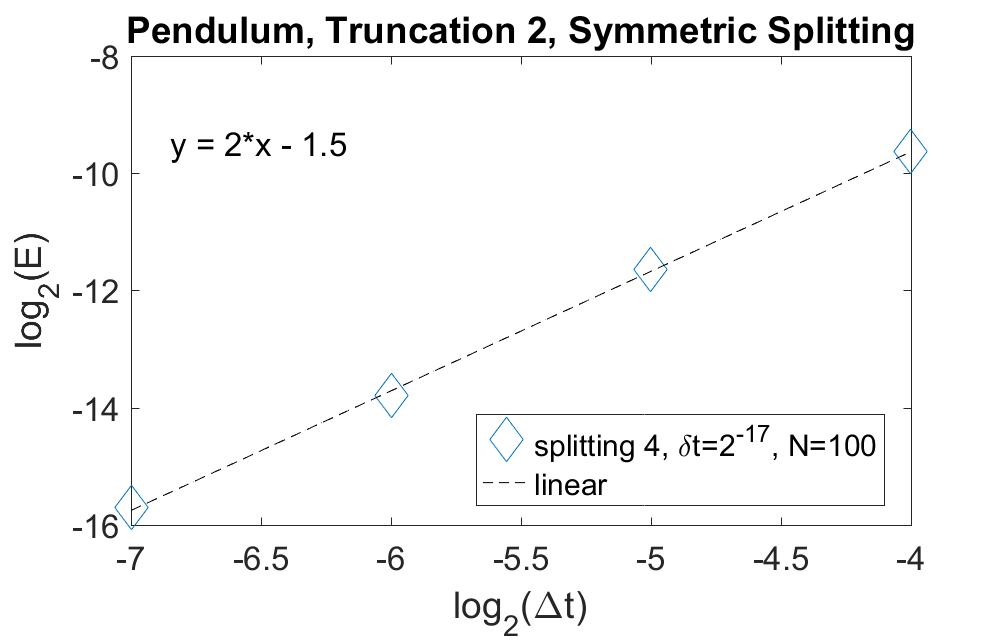}
\includegraphics[scale=.18]{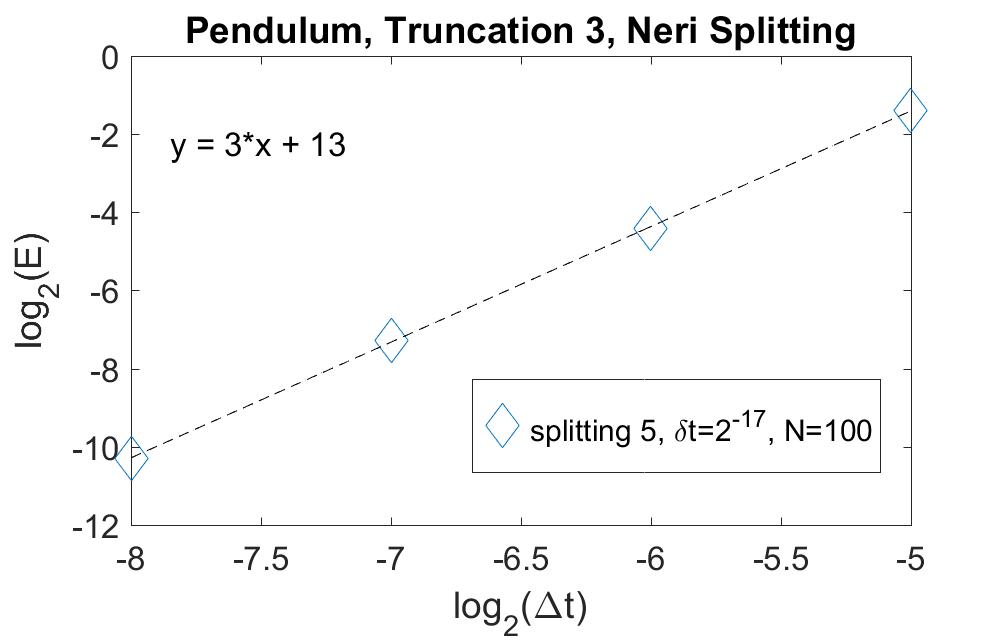}
\caption{The error plot versus the time step for the pendulum model on the log-log scale. From top to bottom: Truncations I, II and III. The step sizes are $\Delta t = 2^{-4},2^{-5},\dots,2^{-8}$. \blue{The number $N$ denotes the number of realizations. Only four values of $\Delta{t}$ appear. Since we are only concerned with asymptotic behavior, we took $\Delta{t}$ as small as possible to obtain the theoretically correct convergence rates. For certain schemes, the values of $\Delta{t}$ needed to be made smaller for this to occur.}}
\label{fig: pend}
\end{center}
\end{figure}

To examine the strong order, we again generate the `exact' solution using the Euler-Maruyama method with very small step size $\delta t= 2^{-17}.$  In order to verify strong convergence, we used 100 realizations. Furthermore, in order to follow the same realization in the implementation of each \blue{algorithm}, we first generate the Brownian motions with small step size, and then the multiple stochastic integrals are evaluated using a numerical quadrature. Notice that this is only for the purpose of examining the strong order of accuracy. In practice, one can sample the integrals using the covariance matrix \eqref{eq: covtrunc3}.
As can be seen from \cref{fig: pend}, the order of accuracy is as expected. \green{The list of numerical tests include
\begin{itemize}
    \item Splitting 1: truncation \rm I with non-symmetric splitting
    \item Splitting 2: truncation \rm I with symmetric splitting
    \item Splitting 3: truncation \rm II with non-symmetric splitting
    \item Splitting 4: truncation \rm II with symmetric splitting
    \item Splitting 5: truncation \rm III with Neri 4th order splitting
\end{itemize}}

\blue{We note that, while Truncation \green{ \rm III } with the Neri splitting has order 3, \green{in this regime} the error is larger than that of the method with order 2. We can see this from the $y$-intercepts on the corresponding graphs. This implies that the 3rd order method has a large prefactor.}

\subsection{A Lennard-Jones cluster}

\begin{figure}[htbp]
\begin{center}
\includegraphics[scale=.17]{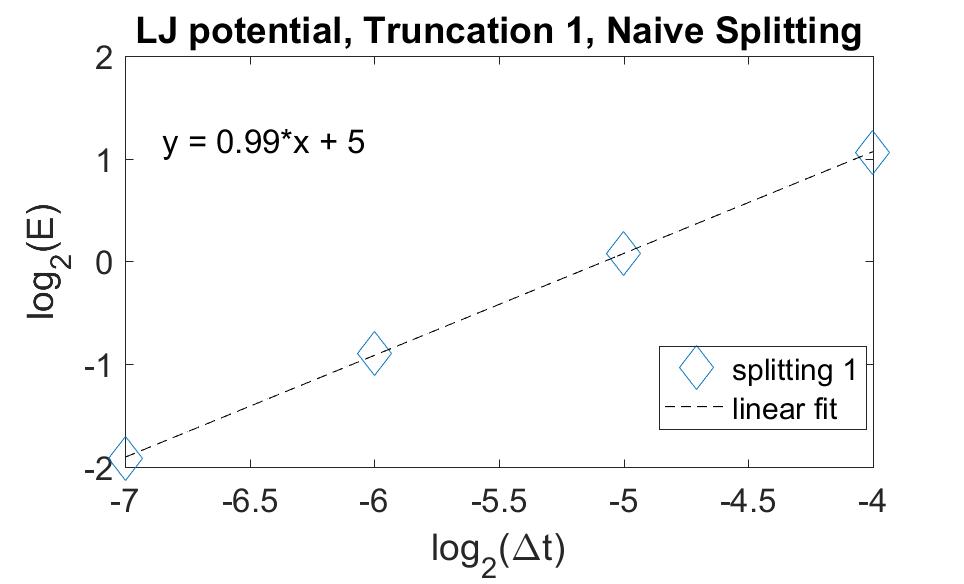}
\includegraphics[scale=.17]{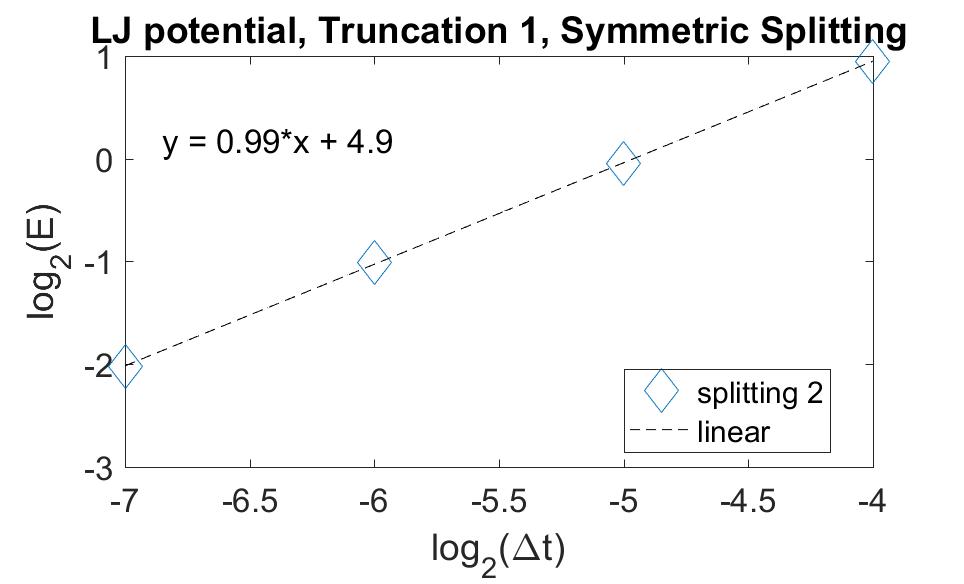}
\includegraphics[scale=.164]{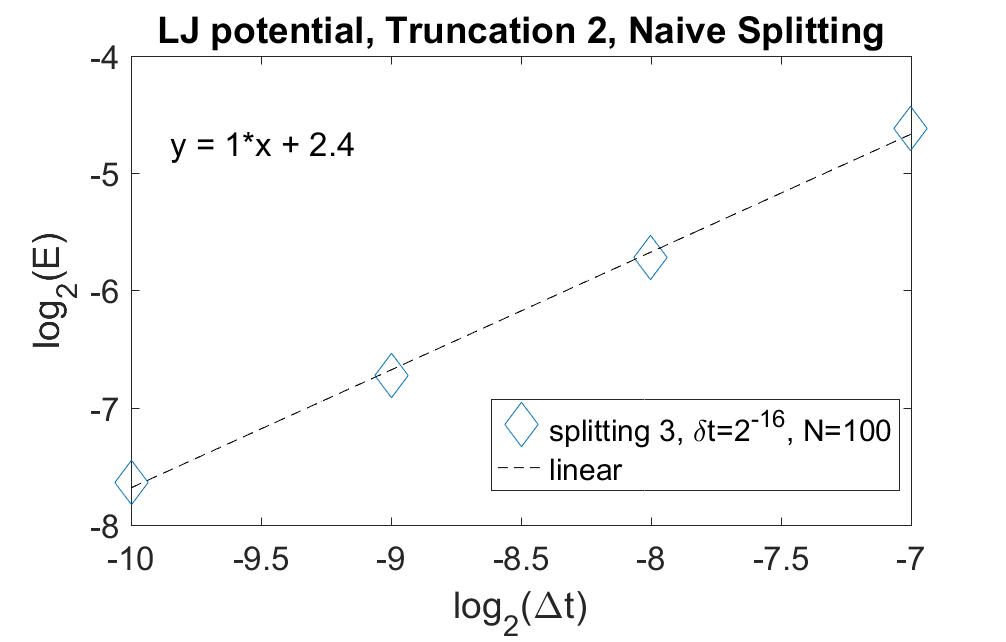}
\includegraphics[scale=.164]{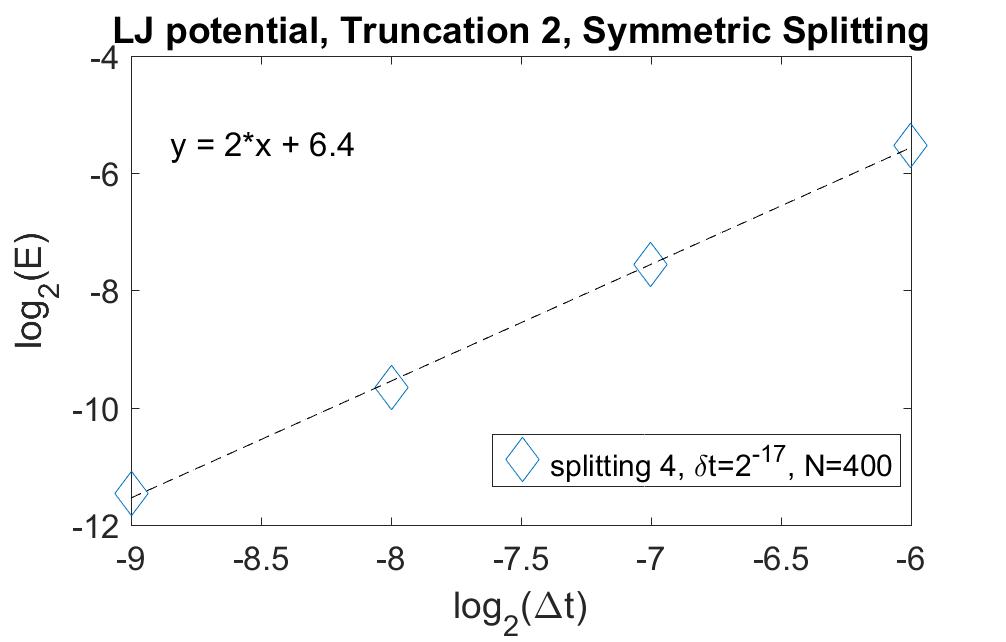}
\includegraphics[scale=.17]{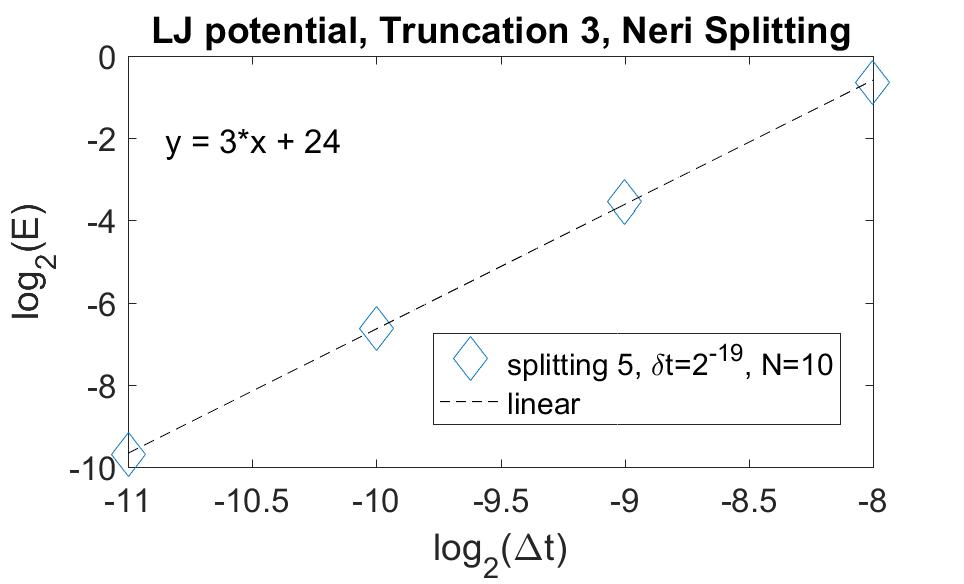}
\includegraphics[scale=.20]{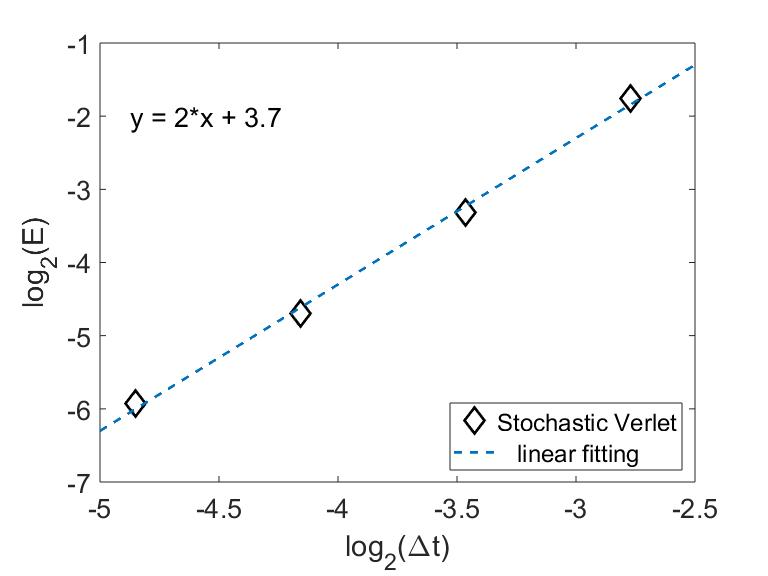}
\caption{The error plot for the LJ-7 cluster on a log-log scale. From top to bottom: Truncations I, II and III, \blue{and the Stochastic velocity Verlet (SVV) method}. The step sizes are $\Delta t=2^{-4},2^{-5},\dots,2^{-11}$. \blue{The number $N$ denotes the number of realizations. Since we are only concerned with asymptotic behavior, we took $\Delta{t}$ as small as possible to obtain the theoretically correct convergence rates.For certain schemes, the values of $\Delta{t}$ needed to be made smaller for this to occur.}}
\label{fig: LJ}
\end{center}
\end{figure}

The $i$th component ($i=1,\dots,n$) of the function $f$ in the Lennard-Jones model is given by
\begin{equation}
f_i(x)=\sum_{j\neq i, j=1}^n \left( 12\left(\frac{1}{r_{ij}}\right)^{13} -6\left( \frac{1}{r_{ij}} \right)^7 \right)\frac{\vec{r}_{ij}}{r_{ij}},
\end{equation}
where $\vec{r}_{ij}=\vec{x}_i-\vec{x}_j\in\R^3$ and $r_{ij}=\|\vec{r}_{ij}\|_2$.
For our simulations, we used seven particles $x_1,\dots,x_7\in\R^3$ and use $x$ to denote the vector $(x_1,\dots,x_7)\in\R^{21}$, so that the total dimension $n$ is 21. Similarly, by $v$ we mean the vector $(v_1,\dots,v_7)\in\R^{21}$ of the velocities of each particle in all dimensions. We take $\Gamma=10I_{21\times 21}$, and $\sigma=\sqrt{2k_B T \Gamma}$, where $k_B T=0.3$, and $T$ is the temperature, not to be confused with the final time $T$. Initially, the atoms are arranged at the vertices of a hexagon and its center. The side length corresponds to the minimum of the Lennard-Jones potential, {$2^{1/6}.$}  Notice that for this model, the function $f(x)$ does not have bounded derivatives unless a cut-off is introduced. Nevertheless, as shown in \cref{fig: LJ}, the strong order of accuracy is still consistent with the results of the analysis. \blue{We note that while the fifth method is indeed of order 3, the error is larger than for the method with order 2, as is the case with the pendulum model.}


\newpage

\section{Summary and discussion}

In this paper, we \green{analyzed the strong convergence of some widely implemented schemes, and} developed new operator-splitting schemes based on Kunita's solution representation. In particular, we obtained algorithms with strong order up to order 3.
This type of convergence is important for predicting the transient stage of the stochastic processes.

There are several remaining challenges in simulating algorithms for Langevin-type of equations. First \cite{schlick1999algorithmic}, there might be multiple scales involved in the force term $f(x)$. In this case, a more appropriate splitting \cite{tuckerman1992reversible} is between the fast and slow forces.
Secondly, the damping  and diffusion coefficients can be position-dependent. Such models arise, for instance,  in the dissipative-particle dynamics (DPD) \cite{warren1998dissipative,Leimkuhler2015}. Finally, there are Langevin equations with strong stiffness, e.g., large damping coefficients. In this case, implicit algorithm are needed. These issues will be addressed in separate works.


\section{Appendix: derivation of It\^o-Taylor expansions with strong convergence rates up to 3}\label{sec: IT}

Due to the lengthy calculations in the analysis, we have included some parts of the proofs in the appendix. These details are useful for the analysis.


\subsection{SDE notation}

An autonomous It\^o stochastic differential equation or SDE can be written formally as,
\begin{equation}
    dz(t)=a(z(t))dt+b(z(t))dW_t,\text{\quad}0\leq t\leq T,
\end{equation}
which means that
\begin{equation}
    z(t)=z(0)+\int_0^t a(z(s))ds+\int_0^t b(z(s))dW_s,\text{\quad}0\leq t\leq T,
\end{equation}
where $\int ds$ denotes the usual Riemann integral and $\int dW_s$ denotes the It\^o stochastic integral,
as defined in \cite[Ch.3]{Oksendal1995} for example, and the $W_t$ denotes a standard $m$-dimensional Wiener process
as defined in \cite[Ch.2]{Oksendal1995}. We use super-scripts for columns or entries. That is, $a=(a^1,\dots,a^d)$,
$b=(b^{ij})_{i,j=1}^d$, $b^j$ is the $j$th column of $b$, and $W=(W^1,\dots,W^m)$.
We call $a=a(z):\R^d\rightarrow\R^d$ the \textit{drift coefficient} and $b=b(z):\R^d\rightarrow \R^{d\times m}$
the \textit{diffusion coefficient}. If $b$ is constant we say the noise is \textit{additive}; otherwise we call
the noise \textit{multiplicative}.

\subsection{Langevin dynamics: notation}

We consider the Langevin dynamics model as a decoupled system of $2n$ SDEs,
\begin{equation}
    \begin{cases}
        dx(t) &= v(t)dt\\
        dv(t) &= \left( f(x(t))-\Gamma v(t) \right)dt + \sigma dW(t),
    \end{cases}
\end{equation}
with drift coefficient
\begin{equation}
    a = \left(\begin{matrix} v \\ f - \Gamma v \end{matrix}\right)\in\R^{2n}
\end{equation}
and constant diffusion coefficient
\begin{equation}
    b = \left(\begin{matrix} 0 \\ \sigma \end{matrix}\right)\in\R^{2n\times n}.
\end{equation}
Here $x,v\in\R^n$ denote position and velocity vectors, $\Gamma,\sigma\in\R^{n\times n}$ are
square constant matrices denoting the friction and diffusion of the system, and
$W\in\R^n$ is the standard Wiener process. Since $b$ is constant, the noise is additive.
The function $f=f(x):\R^n\rightarrow\R^n$ represents a conservative force depending only on position,
which is generally non-linear, and we assume that all second order derivatives of $f$ exist.

\subsection{It\^o-Taylor expansion: notation}

The It\^o-Taylor expansion, see \cite[Ch.5 and 10]{Kloeden1991}, can be viewed as a generalization
of the deterministic Taylor expansion, and of the It\^o formula.
There are weak and strong types of the expansions but we will focus only on strong type in this paper
because we are interested in strong convergence.
We fix a uniform discretization $0=t_0<t_1<\dots<t_N=T$
of the time interval $[0,T]$ with steps $t_k=k\dt$ and uniform step size $\dt=T/N=t_{k+1}-t_k$ for $k=0,1,\dots,N-1$.

Then, the \textit{strong It\^o-Taylor expansion} of order $\gamma\in\{0.5,1.0,1.5,2.0,2.5,3.0\}$ is an explicit one-step method given by $y_0=z(t_0)$
and for $k=1,2,\dots,N-1$,
\begin{equation}
    y_{k+1}=y_k + \sum_{\alpha\in\A_\gamma}g_\alpha(y_k)I_\alpha.
\end{equation}
Here, $\A_\gamma$ is a set of multi-indices $\alpha$, $g_\alpha$ are corresponding coefficient functions of $y_k$,
and $I_\alpha$ are corresponding It\^o stochastic (defined also in \cite[Ch.5]{Kloeden1991}). In any case,
we include the definitions below for completeness so that the reader can understand how we derive
the expansions for the Langevin dynamics model.

Let $\mathcal{M}$ denote the set of all multi-indices $\alpha=(j_1,\dots,j_l)$ with entries $j_1,\dots,j_l\in\{0,1,\dots,m\}$,
length $l(\alpha):=l$, and number of zero entries $n(\alpha)$. Then $\A_\gamma$ is the subset given by,
\begin{equation}
    \A_\gamma:= \{ \alpha\in\mathcal{M}:l(\alpha)+n(\alpha)\leq 2\gamma\text{ or }l(\alpha)+n(\alpha)=\gamma+0.5 \}.
\end{equation}
We will derive $\A_\gamma$ for increasing $\gamma=0.5,1.0,\dots,3.0$,
using the simple fact that $\A_\gamma\subset\A_{\gamma+0.5}$.

Turning now to the coefficient functions $g_\alpha$, we first define differential operators $L^0$ and $L^j$
for $j=1,\dots,m$ as follows:
\begin{equation}
    \begin{aligned}
        L^0&= \sum_{k=1}^n a^k\frac{\partial}{\partial z_k}+\frac{1}{2}\sum_{k,l=1}^n\sum_{j=1}^n
                b^{kj}b^{lj}\frac{\partial^2}{\partial z^k \partial z^l}\\
        L^j&= \sum_{k=1}^n b^{kj}\frac{\partial}{\partial z_k} \text{\quad}j=1,\dots,m.
    \end{aligned}
\end{equation}
The coefficient function $g_\alpha$ corresponding to $\alpha=(j_1,\dots,j_l)\in\mathcal{M}$ is given by,
\begin{equation}
    g_\alpha(z):=L^{j_1}L^{j_2}\dots L^{j_{l-1}}b^{j_l}
\end{equation}
where $b^0:=a$. For the Langevin dynamics model we consider, these operators can be rewritten as,
\begin{equation}
    \begin{aligned}
        L^0&= \sum_{k=1}^n v^k\frac{\partial}{\partial x^k}
                + \sum_{k=1}^n \{f(x)-\Gamma v\}^k\frac{\partial}{\partial v^k}
                +\frac{1}{2}\sum_{k,l=1}^n\sum_{j=1}^n
                        \sigma^{kj}\sigma^{lj}\frac{\partial^2}{\partial v^k\partial v^l}\\
        L^j&= \sum_{k=1}^n \sigma^{kj}\frac{\partial}{\partial v^k}.
    \end{aligned}
\end{equation}
Finally, the stochastic integrals $I_\alpha$ are given by,
\begin{equation}
    I_\alpha:= \int_0^{\dt}\int_0^{s_l}\dots\int_0^{s_3}\int_0^{s_2}du(s_1)du(s_2)\dots du(s_{l-1})du(s_l)
\end{equation}
where each $du(s_i)$ is given by,
\begin{equation}
    du(s_i)= \begin{cases}
                    ds_i &\text{ if }j_i=0 \\
                    dW_{s_i}^{j_i} &\text{ if }j_i\in\{1,\dots,m\}.
    \end{cases}
\end{equation}

\subsection{It\^o-Taylor expansions with strong orders $\gamma = 0.5, 1.0, 1.5, 2.0, 2.5,$ and 3}

First we compute the multi-index sets $\A_\gamma$. For $\gamma=0.5$, we make the simple calculation,
\begin{equation}
    \begin{aligned}
    \A_{0.5}&=\{\alpha:l(\alpha)+n(\alpha)\leq 1\text{ or }l(\alpha)=n(\alpha)=1 \}\\
            &=\{(0),(j):j=1,\dots,m\}.
    \end{aligned}
\end{equation}
The corresponding coefficient functions are,
\begin{equation}
    g_{(0)}=b^0=a=\left(\begin{matrix} v\\ f-\gamma v\end{matrix}\right)
\end{equation}
and
\begin{equation}
    g_{(j)}=b^j=\left(\begin{matrix}0 \\ \sigma^j \end{matrix}\right),\text{\quad}j=1,\dots,m.
\end{equation}
The corresponding stochastic integrals $I_\alpha$ are,
\begin{equation}
    I_{(0)}=\int_0^{\dt}ds_1 = \dt\text{\quad and\quad}I_{(j)}=\int_0^{\dt}dW_{s_1}^j=W_{\dt}^j=:\dw^j.
\end{equation}
Thus, the strong order $\gamma=0.5$ It\^o-Taylor expansion is,
\begin{equation}
    y_{k+1}=y_k+\left(\begin{matrix} v\\ f-\gamma v\end{matrix}\right)\dt
                +\sum_{j=1}^m\left(\begin{matrix}0 \\ \sigma^j \end{matrix}\right)\dw^j,
\end{equation}
which is also called the \textit{Euler} or \textit{Euler-Maruyama} method. As usual $\sigma^j\in\R^d$
denotes the $j$th column of $\sigma\in\R^{d\times m}$, evaluated at $x_k\in\R^d$.

For $\gamma=1.0$, we calculate
\begin{equation}
    \A_1=\{ \alpha:l(\alpha)+n(\alpha)\leq 2 \},
\end{equation}
observing that the equation $l(\alpha)=n(\alpha)=1.5$ is never satisfied. We have
\begin{equation}
    \A_1=\{(0),(j_1),(j_1,j_2):j_1,j_2=1,\dots,m\},
\end{equation}
which is $\A_{0.5}$ plus $(j_1,j_2)$ for $j_1,j_2=1,\dots,m$. The new coefficient functions $g_{(j_1,j_2)}$ are identically zero.
Indeed,
\begin{equation}
    g_{(j_1,j_2)}=L^{j_1}b^{j_2}=\left(\begin{matrix} L^{j_1}(0)\\ L^{j_1}(\sigma^{j_2}) \end{matrix}\right)
            =0\in\R^{2d},
\end{equation}
due to the additive noise.
Therefore, the strong order $\gamma=1.0$ It\^o-Taylor expansion is the same as the one with strong order $\gamma=0.5$.

Consider $\gamma=1.5$. By definition,
\begin{equation}
    \A_{1.5}=\{ \alpha:l(\alpha)+n(\alpha)\leq 3\text{ or }l(\alpha)=n(\alpha)=2 \}.
\end{equation}
Since $\A_1\subset\A_{1.5}$, we just need to find those $\alpha$ for which either
$l(\alpha)+n(\alpha)=3$ or $l(\alpha)=n(\alpha)=2$. We list them now:
\begin{equation}
    (j_1,j_2,j_3),\text{\quad}(0,j_1),\text{\quad}(j_1,0),\text{\quad and\quad}(0,0).
\end{equation}
Fortunately $g_{(j_1,j_2,j_3)}=0$ for $j_1,j_2,j_3=1,\dots,m$, because $g_{(j_1,j_2,j_3)}=L^{j_1}g_{(j_2,j_3)}$ and $g_{(j_2,j_3)}=0$.
Thus the multi-indices $(j_1,j_2,j_3)$ make no contribution to the expansion.
In addition, due to the additive noise the coefficient functions $g_{(0,j)}$ are identically zero as well. Indeed,
\begin{equation}
     g_{(0,j)}= L^0 b^j= 0
\end{equation}
since $b^j=(0,\sigma^j)$ is constant.
On the other hand, we calculate the non-trivial coefficient functions:
\begin{equation}
    \begin{aligned}
    g_{(j,0)}&= L^j a = \left(\begin{matrix} L^j v\\ L^j f-\gamma v \end{matrix}\right)
                    =\left(\begin{matrix} \sigma^j \\ -\gamma\sigma^j \end{matrix}\right),\\
    g_{(0,0)}&= L^0 a = \left(\begin{matrix} L^0 v\\ L^0 \{f-\gamma v\} \end{matrix}\right)
            =\left(\begin{matrix} f-\gamma v\\ \nabla_x(f)v -\gamma \{f-\gamma v\} \end{matrix}\right).
    \end{aligned}
\end{equation}
The corresponding It\^o integrals are:
\begin{align}
    I_{(j,0)}&=\int_0^{\dt}\int_0^{s_2}dW_{s_1}^j ds = \int_0^{\dt}W_s ds\\
    I_{(0,0)}&=\int_0^{\dt}\int_0^{s_2}ds_1 ds_2 = \frac{\dt^2}{2}.
\end{align}
Thus, the strong order $\gamma=1.5$ It\^o-Taylor expansion is,
\begin{align}
    y_{k+1}&=y_k+\left(\begin{matrix} v\\ f-\gamma v\end{matrix}\right)\dt
                +\sum_{j=1}^m\left(\begin{matrix}0 \\ \sigma^j \end{matrix}\right)\dw^j \\
                &+\sum_{j=1}^m\left(\begin{matrix} \sigma^j \\ -\gamma\sigma^j \end{matrix}\right)\int_0^{\dt}W_s^j ds\\
                &+\left(\begin{matrix} f-\gamma v\\ \nabla_x(f)v -\gamma \{f-\gamma v\} \end{matrix}\right)\frac{\dt^2}{2}.
\end{align}
The notation $\nabla_x(f)$ denotes the $d\times d$ Jacobian matrix of $f$ evaluated at $x_k$, with derivatives taken with respect to $x$.

We move on to $\gamma=2.0$. The hierarchical set $\A_2$ is,
\begin{equation}
    \A_2=\{ \alpha:l(\alpha)+n(\alpha)\leq 4 \},
\end{equation}
noting that the equation $l(\alpha)=n(\alpha)=2.5$ is never satisfied.
We also recall that $\A_{1.5}\subset\A_2$, and so we need only look for new multi-indices $\alpha$
which satisfy $l(\alpha)+n(\alpha)=4$, which are,
\begin{equation}
    (j_1,j_2,0),\text{\quad}(j_1,0,j_2),\text{\quad}(0,j_1,j_2),\text{\quad and\quad}(j_1,j_2,j_3),\text{\quad for }j_1,j_2,j_3=1,\dots,m.
\end{equation}
The coefficient functions $g_{(0,j_1,j_2)}$ and $g_{(j_1,j_2,j_3)}$ are identically zero,
since $g_{(j_1,j_2)}$ is zero. In addition, $g_{(j_1,j_2,0)}=0$ because $g_{(j,0)}=(\sigma^j, -\gamma\sigma^j)$ is constant.
That leaves us only with $g_{(j_1,0,j_2)}$, which is also zero because $g_{(0,j_2)}=0$.
\begin{equation}
    g_{(j_1,0,j_2)}=L^{j_1}g_{(0,j_2)}=\left(\begin{matrix} L^{j_1}0\\ L^{j_1}\nabla_x(\sigma^{j_2})v \end{matrix}\right)
        =\left(\begin{matrix} 0\\ \nabla_x(\sigma^{j_2})\sigma^{j_1} \end{matrix}\right).
\end{equation}
Therefore there are no new non-trivial coefficient functions, and thus the strong order $\gamma=2.0$ It\^o-Taylor
expansion is the same as the one with order $\gamma=1.5$.

As for $\gamma=2.5$, we calculate
\begin{equation}
\begin{aligned}
    \A_{2.5}-\A_2&= \{ \alpha\in\A_{2.5}:\alpha\notin\A_2 \}\\
        &= \{ (0,0,j_1),(0,j_1,0),(j_1,0,0),(0,0,0)\\
            &(j_1,j_2,j_3,0),(j_1,j_2,0,j_3),(j_1,0,j_2,j_3),(0,j_1,j_2,j_3)\\
            &(j_1,j_2,j_3,j_4,j_5):j_i=1,\dots,m\text{ and }i=1,\dots,5 \}.
\end{aligned}
\end{equation}
Next we find the vanishing coefficient functions $g_\alpha=0$ for $\alpha\in\A_{2.5}-\A_2$ in the case $m=1$,
which generalizes naturally to $m>1$. Since $g_{(1)}=b=(0,\sigma)$ is constant,
the coefficient functions $g_{(0,0,1)},g_{(1,1,0,1)},g_{(0,1,1,1)},g_{(1,1,1,1,1)}$ vanish:
\begin{align}
    g_{(0,0,1)}& = L^0L^0 g_{(1)}=0\\
    g_{(1,1,0,1)}&= L^1L^1L^0 g_{(1)}=0\\
    g_{(1,0,1,1)}&= L^1 L^0 L^1 g_{(1)}=0\\
    g_{(0,1,1,1)}&= L^0 L^1 L^1 g_{(1)}=0\\
    g_{(1,1,1,1,1)}&= L^1 L^1 L^1 L^1 g_{(1)}=0.
\end{align}
What about the indices $\alpha=(1,1,1,0),(0,1,0),(1,0,0)$ and $(0,0,0)$?
Well, we have already seen that $g_{(1,0)}=(\sigma,-\gamma\sigma)$ is constant.
Therefore, $g_{(1,1,1,0)}=g_{(0,1,0)}=0$, since
\begin{align}
    g_{(1,1,1,0)}&= L^1 L^1 g_{(1,0)}=0\\
    g_{(0,1,0)}&= L^0 g_{(1,0)}=0.
\end{align}
We have two indices remaining: $\alpha=(1,0,0)$ and $(0,0,0)$, and the corresponding coefficient functions
do not vanish, as we will see. We recall that,
\begin{align}
    g_{(0,0)}&= \left(\begin{matrix} f-\gamma v\\ f'v-\gamma\{f-\gamma v\} \end{matrix}\right)
\end{align}
for additive noise, and the differential operators are,
\begin{align}
    L^0&= v\partial_x +\{f-\gamma v\}\partial_v\\
    L^1&= \sigma\partial_v.
\end{align}
Therefore,
\begin{align}
    g_{(1,0,0)}&= L^1 g_{(0,0)}\\
        &= \left(\begin{matrix} \sigma\partial_v \{f-\gamma v\} \\
                \sigma\partial_v\left( f'v-\gamma \{f-\gamma v\} \right) \end{matrix}\right)\\
        &= \left(\begin{matrix} -\gamma\sigma \\ f' \sigma +\gamma^2 \sigma \end{matrix}\right),
\end{align}
and
\begin{align}
    g_{(0,0,0)}&= L^0 g_{(0,0)}\\
        &=\left(\begin{matrix} v\partial_x (f-\gamma v)+\{f-\gamma v\}\partial_v(f-\gamma v)\\
            v\partial_x[f'v-\gamma\{f-\gamma v\}]+\{f-\gamma v\}\partial_v[f'v-\gamma\{f-\gamma v\}] \end{matrix}\right)\\
        &= \left(\begin{matrix} f'v-\gamma\{f-\gamma v\}\\
                f''v^2-\gamma f'v + f'\{f-\gamma v\}+\sigma^2\{f-\gamma v\} \end{matrix}\right).
\end{align}
The corresponding functions for the case $m>1$ are,
\begin{align}
    g_{(j,0,0)}&= \left(\begin{matrix} -\gamma\sigma^j \\ \nabla_x(f)\sigma^j +\gamma^2 \sigma^j \end{matrix}\right) \text{\quad and}\\
    g_{(0,0,0)}&= \left(\begin{matrix}\nabla_x(f)v-\gamma\{f-\gamma v\} \\ \nabla_x[\nabla_x(f)v]v-\gamma \nabla_x(f)v+\nabla_x(f)\{f-\gamma v\}+\sigma^2\{f-\gamma v\} \end{matrix}\right),
\end{align}
which is a straightforward, but mildly unpleasant calculation. The related stochastic It\^o integrals are,
\begin{align}
    I_{(1,0,0)}&= \int_0^{\dt}\int_0^{s_3}\int_0^{s_2}dW_{s_1}ds_2 ds_3\text{\quad and}\\
    I_{(0,0,0)}&= \int_0^{\dt}\int_0^{s_3}\int_0^{s_2}dW_{s_1}ds_2 ds_3=\frac{\dt^3}{3!}.
\end{align}
Now we can add these terms to the strong order $\gamma=2$ expansion to get the one with strong order $\gamma=2.5$:
\begin{align}
    y_{k+1}&=y_k+\left(\begin{matrix} v\\ f-\gamma v\end{matrix}\right)\dt
                +\sum_{j=1}^m\left(\begin{matrix}0 \\ \sigma^j \end{matrix}\right)\dw^j \\
                &+\sum_{j=1}^m\left(\begin{matrix} \sigma^j \\ -\gamma\sigma^j \end{matrix}\right)\int_0^{\dt}W_s^j ds\\
                &+\left(\begin{matrix} f-\gamma v\\ \nabla_x(f)v -\gamma \{f-\gamma v\} \end{matrix}\right)\frac{\dt^2}{2}\\
                &+\sum_{j=1}^m\left(\begin{matrix} -\gamma\sigma^j \\ \nabla_x(f)\sigma^j +\gamma^2 \sigma^j \end{matrix}\right)\int_0^{\dt}\int_0^{s_3}\int_0^{s_2}dW_{s_1}^jds_2 ds_3\\
                &+\left(\begin{matrix}\nabla_x(f)v-\gamma\{f-\gamma v\} \\ \nabla_x[\nabla_x(f)v]v-\gamma \nabla_x(f)v+\nabla_x(f)\{f-\gamma v\}+\sigma^2\{f-\gamma v\} \end{matrix}\right)\frac{\dt^3}{3!}.
\end{align}
It turns out that the strong order 3 and strong order 2.5 methods are the same for additive noise, and here is why.
First, we have the hierarchical set,
\begin{equation}
    \A_3=\{\alpha\in\mathcal{M}: l(\alpha)+n(\alpha)\leq 6 \},
\end{equation}
noting that the condition $l(\alpha)=n(\alpha)=3.5$ is not possible. We focus on the case $m=1$, that is,
the Wiener process is one dimensional. We observe that the new multi-indices are,
\begin{align}
    \A_3-\A_{2.5}&= \{\alpha\in\mathcal{M}: l(\alpha)+n(\alpha)=6\text{ and }l(\alpha)=n(\alpha)=3\text{ is false} \}\\
        &=\{ (0,0,1,1),(0,1,0,1),(0,1,1,0),(1,0,0,1),(1,0,1,0),(1,1,0,0)\\
            &(1,1,1,1,0),(1,1,1,0,1),(1,1,0,1,1),(1,0,1,1,1),(0,1,1,1,1)\\
            &(1,1,1,1,1,1) \}.
\end{align}
Since $g_{(1)}$ and $g_{(1,0)}$ are constant, it is clear that the coefficient functions vanish
for every $\alpha$ in $\A_3-\A_{2.5}$ except for possibly $(1,1,0,0)$.
But $g_{(1,1,0,0)}=L^1 g_{(1,0,0)}$, and $g_{(1,0,0)}$ is constant with respect to $v$, and $L^1$ differentiates with respect to $v$.
Therefore $g_{(1,1,0,0)}=0$ as well.

In conclusion, every new $\alpha\in\A_3-\A_{2.5}$ corresponds to a vanishing coefficient function $g_\alpha$,
and hence the strong order 3 method is the same -- for additive noise -- as the strong order 2.5 method.
This also concludes our analysis of strong It\^o-Taylor expansions.


\bibliographystyle{siamplain}
\bibliography{Num4Lgv,nemd}
\end{document}